\documentclass[11pt]{amsart}
\usepackage{graphicx}
\usepackage{amssymb,amscd,amsthm,amsxtra}
\usepackage{latexsym}
\usepackage{epsfig}
\usepackage{amsmath}
\usepackage{bm}
\usepackage[usenames,dvipsnames,svgnames,table]{xcolor}
\usepackage[linktocpage=true,colorlinks=true,linkcolor=Blue,citecolor=BrickRed,urlcolor=RoyalBlue]{hyperref}
\usepackage[alphabetic]{amsrefs}
\usepackage[colorinlistoftodos,prependcaption,textsize=tiny]{todonotes}
\usepackage[mathscr]{euscript} 
\usepackage{mathrsfs} 
\usepackage{pxfonts}
\usepackage[varg]{txfonts}
\usepackage{esint}
\usepackage{tikz}
\usepackage{verbatim}
\usepackage{setspace}
\usepackage[normalem]{ulem}

\numberwithin{equation}{section}

\title[Sharp discretization of $L_p$ norms]{\texorpdfstring{Sharp discretization of $L_p$ norms}{Sharp discretization of Lp norms}}
\author[D. M. Pellegrino]{Daniel M. Pellegrino}
\address{Departamento de Matem\'atica, Universidade Federal da Para\'iba,
Campus I, 58051-900 Jo\~ao Pessoa, PB, Brazil}
\email{daniel.pellegrino@academico.ufpb.br}
\author[E. V. Teixeira]{Eduardo V. Teixeira}
\address{Department of Mathematics, Oklahoma State University, 401 Mathematical
Sciences, Stillwater, OK 74078, USA}
\email{eduardo.teixeira@okstate.edu}

\newtheorem{theorem}{Theorem}
\newtheorem{lemma}{Lemma}[section]
\newtheorem{corollary}{Corollary}[section]
\newtheorem{proposition}{Proposition}[section]

\newcommand{\R}{\mathbb{R}}
\newcommand{\E}{\mathbb{E}}
\newcommand{\PP}{\mathbb{P}}

\newcommand{\eps}{\varepsilon}

\newcommand{\intav}[1]{\mathchoice 
  {\mathop{\vrule width 6pt height 3 pt depth -2.5pt \kern -8pt \intop}\nolimits_{\kern -6pt#1}} 
  {\mathop{\vrule width 5pt height 3 pt depth -2.6pt \kern -6pt \intop}\nolimits_{#1}}
  {\mathop{\vrule width 5pt height 3 pt depth -2.6pt \kern -6pt \intop}\nolimits_{#1}}
  {\mathop{\vrule width 5pt height 3 pt depth -2.6pt \kern -6pt \intop}\nolimits_{#1}}}

\numberwithin{equation}{section}
\newtheorem{mainresult}{Theorem}

\newcommand{\Sph}{\mathbb S}
\DeclareMathOperator{\supp}{supp}
\newcommand{\op}{\mathrm{op}}
\newcommand{\eqid}[1]{\leqno\hbox{\refstepcounter{equation}\label{#1}\normalfont(\theequation)}}
\hypersetup{pdftitle={Sharp discretization of Lp norms},
  pdfauthor={Daniel M. Pellegrino and Eduardo V. Teixeira},bookmarksnumbered=true}

\begin{document}

\subjclass[2020]{Primary 46B28, 46B07. Secondary 46B09, 46E30, 60B20}
\keywords{Absolutely summing operators, discretization of power moments,
random Euclidean sections, partial coloring}
\date{}

\begin{abstract}
We remove the logarithmic factor in the embedding theorem of Bourgain, Lindenstrauss and Milman for finite-dimensional subspaces of $L_p$ and settle the question posed by Johnson and Schechtman in 1994 on the approximation of $p$-summing norms. Both results are special cases of our main theorem, which establishes a discretization of $L_p$ norms for every $p>2$, with optimal dimensional order and exact preservation of prescribed linear constraints.
\end{abstract}

\maketitle

\section{Introduction}

Given $2<p<\infty$ and $0<\varepsilon<1/2$, a natural problem is to determine the smallest integer $m$ such that every $n$-dimensional subspace of real $L_p(0,1)$ embeds linearly into $\ell_p^m$ with distortion at most $1+\varepsilon$. In 1989, Bourgain, Lindenstrauss and Milman \cite{BLM}*{Theorem~7.3} showed that
\begin{equation}\label{intro:BLM}
 m\le C_{p,\varepsilon}n^{p/2}\log(2+n).
\end{equation}
The exponent $p/2$ is optimal. The question is whether the logarithmic factor in \eqref{intro:BLM} is necessary.

The same dimensional problem arises in the theory of absolutely summing operators. In 1994, Johnson and Schechtman \cite{JS} refined the discretization argument of Bourgain, Lindenstrauss and Milman to approximate the $p$-summing norm of an operator with $n$-dimensional domain. For $p>2$ and fixed $0<\varepsilon<1$, they showed that $C_{p,\varepsilon}n^{p/2}\log^3(2+n)$ vectors suffice for a $(1+\varepsilon)$-approximation of the norm. They proved that the power of $n$ is optimal and explicitly asked in the introduction whether a logarithmic factor is needed. The operator problem requires more than the existence of an almost isometric copy in a finite-dimensional $\ell_p$ space: the discretization must also control the numerator in the quotient defining the $p$-summing norm.

The two problems have subsequently been studied from different directions. In 2003, Zvavitch \cite{Zv}*{Theorem~1.1(ii)} reduced the logarithmic loss in the operator problem, obtaining a $(1+\varepsilon)$-approximation with
\begin{equation}\label{intro:Zv}
 m\le \left\lceil C_p\varepsilon^{-5}n^{p/2}\log(2+\varepsilon^{-5}n)\right\rceil.
\end{equation}
Johnson and Schechtman \cite{JSlocal} studied the geometric problem for particular classes of subspaces of $L_p$. In 2011, Schechtman \cite{Even} removed the logarithm from \eqref{intro:BLM} for even integers $p$, using the quadratic sparsification theorem of Batson, Spielman and Srivastava \cite{BSS}. The argument passes from an $L_p$ norm to a quadratic form on the span of products of $p/2$ functions, a construction that relies on $p$ being even. The significance of this advance in the geometric problem is also discussed in Naor's Bourbaki exposition \cite{Naor}*{Section~5.2}.

We remove the logarithmic factor in \eqref{intro:BLM} for every real exponent $p>2$ and settle the question of Johnson and Schechtman in the same range. More precisely, for each $p>2$, $\sigma>0$ and $0<\varepsilon<1$, the integer
\begin{equation}\label{intro:optimal-m}
 m=\left\lceil C_{p,\sigma}\varepsilon^{-2-\sigma}n^{p/2}\right\rceil
\end{equation}
suffices both for a linear embedding of any $n$-dimensional subspace of $L_p$ into $\ell_p^m$ with distortion at most $1+\varepsilon$ (Corollary~\ref{full:embedding}) and for a $(1+\varepsilon)$-approximation of the $p$-summing norm of any operator with $n$-dimensional domain (Corollary~\ref{full:witness}). Both assertions hold over the real and complex fields. The order in $n$ is optimal for fixed $p$ and $\varepsilon$; we do not determine the optimal dependence on $\varepsilon$.

Both estimates follow from the same discretization theorem for power moments. Let $a_1,\ldots,a_N$ be linear forms on a finite-dimensional
space $V$. For nonnegative weights $w=(w_1,\ldots,w_N)$, write
$$
 F(y)=\sum_{i=1}^{N}|a_i(y)|^p,\qquad
 F_w(y)=\sum_{i=1}^{N}w_i|a_i(y)|^p .
$$
We seek a sparse $w$ for which $F_w/F$ is uniformly close to one wherever
$F>0$. Each index may also carry a vector $c_i$ of prescribed costs.
The forms may repeat or vanish, and distinct indices carrying the same
form may have different costs.

Our main theorem encompasses both problems. The case $s=0$ yields the optimal dimensional bound in the embedding problem of Bourgain, Lindenstrauss and Milman, while the case $s=1$, with $c_i=\|Tx_i\|^p$, settles the question of Johnson and Schechtman. For arbitrary $s$, prescribed linear constraints are preserved exactly, at an additive cost in the support size.

\begin{mainresult}[Discretization with prescribed linear costs]\label{full:rows}
Fix $2<p<\infty$, $\sigma>0$ and $0<\eta<1$. Suppose $\dim V=n\ge1$ over
$\mathbb R$ or $\mathbb C$. Given $a_1,\ldots,a_N\in V^*$ and
$c_1,\ldots,c_N\in\mathbb R^s$, with $s\ge0$, there exist nonnegative
weights $w_1,\ldots,w_N$ such that
$$
 \begin{aligned}
 |\supp w|&\le C_{p,\sigma}\eta^{-2-\sigma}n^{p/2}+128(s+1),\\
 \sum_iw_i c_i&=\sum_i c_i,\\
 (1-\eta)F(y)&\le F_w(y)\le(1+\eta)F(y)\qquad (y\in V).
 \end{aligned}
$$
The constant depends only on $p$ and $\sigma$.
\end{mainresult}

For a bounded linear map $T:X\to Y$ and a positive integer $m$, let
$$
 \pi_p^{(m)}(T)=
 \sup_{(x_1,\ldots,x_m)\ne0}
 \frac{\big(\sum_{i=1}^m\|Tx_i\|^p\big)^{1/p}}
 {\displaystyle\sup_{\varphi\in B_{X^*}}
       \big(\sum_{i=1}^m|\varphi(x_i)|^p\big)^{1/p}},
 \qquad \pi_p(T)=\sup_{m\ge1}\pi_p^{(m)}(T).
 \eqid{intro:pi}
$$
Here $B_{X^*}$ denotes the closed dual unit ball. A family of fewer than $m$
vectors can be completed by zeros. The numerator in \eqref{intro:pi} is
the reason for allowing a cost in Theorem~\ref{full:rows}.

\begin{corollary}[Finite witnesses]\label{full:witness}
Let $2<p<\infty$, $\sigma>0$ and $0<\varepsilon<1$. For every bounded linear operator $T:X\to Y$
between real or complex Banach spaces with $\dim X=n\ge1$,
$$
 \pi_p(T)\le(1+\varepsilon)\pi_p^{(m)}(T),
 \qquad m=\left\lceil C_{p,\sigma}\varepsilon^{-2-\sigma}n^{p/2}\right\rceil .
$$
\end{corollary}

\begin{proof}
For a family $(x_i)$ with positive numerator, apply Theorem~\ref{full:rows} on $X^*$ to
$a_i(\varphi)=\varphi(x_i)$ and $c_i=\|Tx_i\|^p$. The vectors $z_i=w_i^{1/p}x_i$ satisfy
$$
 \sum_i\|Tz_i\|^p=\sum_i\|Tx_i\|^p,\qquad
 \sup_{\varphi\in B_{X^*}}\sum_i|\varphi(z_i)|^p
 \le(1+\varepsilon)\sup_{\varphi\in B_{X^*}}\sum_i|\varphi(x_i)|^p .
$$
The original quotient is therefore at most $(1+\varepsilon)^{1/p}\pi_p^{(m)}(T)$. The additive
constant in the support bound is absorbed in $C_{p,\sigma}$, uniformly for $n\ge1$ and
$0<\varepsilon<1$. Taking the supremum proves the assertion, including $T=0$.
\end{proof}

The dimension hypothesis in Corollary~\ref{full:witness} is essential. For $1<p<\infty$,
$p\ne2$, Johnson and Schechtman showed that any bound valid for all operators of rank $n$, with
unrestricted domain dimension, requires more than polynomially many vectors
\cite{JS}*{Section~IV}.

Each additional cost
incurs only an additive term in the support. For comparison, Cohen and
Peng study algorithmic row sampling by Lewis weights \cite{CP}, and Dai,
Kosov and Murasko obtain estimates for independent, equal-weight samples
under a Nikol'skii inequality \cite{DKM}*{Corollary~1.3}.

To see why the power of $n$ is sharp, let $I_n$ denote the identity on
real $\ell_2^n$, and let $c_{n,p}$ be the $p$th absolute moment of the
first coordinate of a uniform unit vector. For any family $(x_j)$,
averaging $\sum_j|\langle u,x_j\rangle|^p$ over $u\in\Sph^{n-1}$
gives $c_{n,p}\sum_j|x_j|^p$. Conversely, finite spherical quadratures
approximate the constant moment $c_{n,p}$ uniformly in the test
direction, by compactness. This gives Gordon's identity
$\pi_p(I_n)^p=c_{n,p}^{-1}\asymp_p n^{p/2}$ \cite{Gordon}.
Testing in the direction of a longest vector gives
$\pi_p^{(m)}(I_n)^p\le m$.

\subsection*{What happens at the critical cardinality}
The order of $m$ leaves open how much of the norm a family of this size
can recover. Write $\mu_p=\mathbb E|g|^p$ for a standard real Gaussian
$g$. If $m\sim\alpha n^{p/2}$, $0<\alpha<\infty$, independent uniform
unit vectors recover, in probability, the proportion
$$
 \left(\frac{\alpha\mu_p}{1+\alpha\mu_p}\right)^{1/p}
$$
of $\pi_p(I_n)$. The term $1$ in the denominator has a geometric
origin. Testing in the direction of one sampled vector makes its own
contribution equal to one, in addition to the background moment
$\alpha\mu_p$. Theorem~\ref{th:fixed} describes the directions on which
this excess is attained. Allowing the $m$ vectors to be selected from a larger sample
and reweighted optimally gives the same limiting proportion whenever
the sample size $M$ satisfies $M\ge m$ and $\log M=o(n)$
(Theorem~\ref{th:adaptive}).

For a lower bound valid in the operator problem, the witnesses must
be free to range over the entire domain. Theorem~\ref{thm:operator}
shows that the diagonal contribution persists for random evaluation
maps into $\ell_\infty^M$, even with this freedom. Consequently every
coefficient $C_{p,\varepsilon}$ valid for all real operators in
Corollary~\ref{full:witness}, with $m=\lceil C_{p,\varepsilon}n^{p/2}\rceil$,
satisfies
$$
 C_{p,\varepsilon}\ge\frac{1}{\mu_p((1+\varepsilon)^p-1)}
       \sim\frac{1}{p\mu_p\varepsilon}\qquad(\varepsilon\downarrow0).
$$
The optimal dependence on $\varepsilon$ for general operators remains
open. Appendix~\ref{sec:refinements} gives an
$O_p(\varepsilon^{-2}\log(2/\varepsilon)n^{p/2})$ bound on Hilbert
domains and a quadratic dependence on inverse accuracy for even
moments with prescribed costs.

The extremal spherical moment also determines the distortion of a
Haar-distributed $k$-dimensional section of $\ell_p^N$. For fixed $p>2$
and $k/N^{2/p}\to\beta\in(0,\infty)$, its Banach--Mazur distance to
Euclidean space converges in probability to
$(1+\beta^{p/2}/\mu_p)^{1/p}$. At logarithmic exponents, the norm in
each Gaussian direction is governed by the larger of a deterministic
Euclidean term and the largest coordinate. Theorems~\ref{th:uniform}
and~\ref{th:main} give the distortion on both sides of $p=2\log N$.
Optimizing over all Euclidean structures on the section yields the same
leading constant as the inherited Euclidean norm.

\subsection*{The discretization argument}
Choose Lewis coordinates and suppose the current family has $N$ rows
of rank $r$. At size $W=N/(8r^{p/2})$, it suffices to remove a fixed
proportion of the rows at a relative cost of order $W^{-1/(2+\sigma)}$.
These errors form a geometric sum when the reductions are read in reverse.
The difficulty lies in controlling every test direction with a single
choice of rows.

We project the rows onto a nested sequence of moment balls. Averaging
the odd and even projections separately gives two families whose
difference has row norms of order $J^{-1}$ after $J$ pairs. This
difference still has a higher moment bound on the new test vectors,
so the construction can be repeated. The shorter rows at each
successive stage require fewer coefficients. To benefit from that
decrease one must count each chosen index using the degree of the
graph at the stage where the index was chosen. Section~\ref{sec:support-count}
does this by truncating the approximations on their supports. A single
partial coloring then controls all dyadic heights and resolutions.
The prescribed costs are orthogonality constraints on the coloring.
The assumption $p>2$ enters in summing the estimates over heights in
Section~\ref{full:halving}. All logarithms are natural unless a base
is indicated.

\begingroup
\value{tocdepth}=1\relax
\tableofcontents
\endgroup

\section{Lewis position and partial coloring}\label{full:normalization}

The proof rests on a reduction that can be iterated. A surviving row absorbs the $p$th root of
its new weight, while its cost absorbs the weight itself. After each reduction we therefore have
another instance of the same problem.

\begin{proposition}[Reduction of the support]\label{det:reduction}
Let a real family have $N$ labels, rank $r\ge1$, and $s$ real linear costs. If $N\ge128s$ and
$W=N/(8r^{p/2})\ge1$, there are nonnegative multipliers, supported on at most $27N/32$ original
labels, that preserve the total cost exactly and change every $p$th moment by at most
$C_{p,\sigma}W^{-1/(2+\sigma)}$ times its original value.
\end{proposition}

Pass to the quotient of the domain by the common kernel of the rows and choose real Euclidean
coordinates on its dual. The rows now span $\mathbb R^r$, with $r\le n$. We treat complex forms at
the end of Section~\ref{full:conclusion}. When $r=0$ the moment vanishes. Conic Carath\'eodory
represents the cost sum by at most $s$ generators. Indeed, given a dependence among generators
with positive weights, subtract a multiple of that dependence until the first weight vanishes.
If $s=0$ all weights may be taken zero. Henceforth $r\ge1$.

An invertible linear change gives Lewis coordinates \cite{Lewis}.
In fact, minimize $\sum_i(a_i^{\top}Ha_i)^{p/2}$ over positive
definite $H$ of determinant one. The rows span, so
$\inf_{|v|=1}\max_i|\langle a_i,v\rangle|>0$. The objective tends
to infinity with the largest eigenvalue, and a minimizer exists.
Differentiation along $H^{1/2}\exp(tS)H^{1/2}$ for traceless symmetric
$S$ shows that, after replacing $a_i$ by $H^{1/2}a_i$,
$$
 \sum_i |a_i|^{p-2}a_ia_i^{\top}=\frac Br I_r,
 \qquad B=\sum_i|a_i|^p .
$$
For $a_i\ne0$, set $u_i=a_i/|a_i|$ and $\nu_i=|a_i|^p/B$. Isotropy gives $\sum_i\nu_i|\langle
u_i,y\rangle|^2=|y|^2/r$. Convexity of $t^{p/2}$ on the nonnegative half-line now gives
$\sum_i|\langle a_i,y\rangle|^p\ge B r^{-p/2}|y|^p$.

Let $b_0=8B/N$ and split $i$ into
$k_i=\max\{1,\lceil |a_i|^p/b_0\rceil\}$ copies. Give each copy
the form $a_i/k_i^{1/p}$ and cost $c_i/k_i$. Their total number $M$
satisfies $N\le M\le9N/8$. Divide their forms by $b_0^{1/p}$ and
write the normalized rows as $a_1,\ldots,a_M$ again. Then
$$
 |a_i|\le1,\qquad F(y)=\sum_{i=1}^M|\langle a_i,y\rangle|^p,
 \qquad W=\frac{N}{8r^{p/2}},\qquad
 K=\{y:F(y)\le W\}\subset B_2^r .
$$
Indeed, $k_i\le1+|a_i|^p/b_0$ and $B/b_0=N/8$. Splitting leaves the moment unchanged. The common
rescaling then gives $F(y)\ge W|y|^p$ by isotropy. The body $K$ is compact, convex and
symmetric, and $F=W$ on its
boundary. A relative estimate there extends to every $y$ by homogeneity.
The reduction will be used when $W\ge1$ and $s\le M/128$. Below, $|\cdot|$ denotes the Euclidean
norm or the absolute value, and $\ell_q$ norms of coordinate arrays use counting measure.

We use the partial coloring theorem of Lovett and Meka
\cite{LM}*{Theorem~4 in the preprint}, in the following form.

\begin{lemma}[Scalar partial coloring]\label{det:LM}
Let $\mathcal V$ be a finite or countable collection of vectors in $\mathbb R^M$. Given
nonnegative thresholds $(\lambda_v)_{v\in\mathcal V}$ with
$\sum_v\exp(-\lambda_v^2/16)\le M/16$, there exists
$\chi\in[-1,1]^M$ such that $|\langle\chi,v\rangle|\le\lambda_v|v|$
for every $v$ and at least $M/2$ coordinates have modulus one.
\end{lemma}

\begin{proof}
Apply the finite partial coloring theorem with starting point zero and with boundary tolerance
tending to zero. A subsequence converges in the cube, and a further subsequence has the same set
of at least $M/2$ nearly saturated coordinates. The scalar inequalities are closed and the
limiting coordinates have modulus one. For countably many constraints, use this conclusion for
their first $k$ members and pass to a convergent subsequence as $k\to\infty$. Every fixed
constraint holds eventually along the subsequence. The set of points with at least $M/2$ saturated
coordinates is a finite union of closed faces of the cube.
\end{proof}

\section{Interlaced projections and the residual moment}
\label{sec:weighted-interlacing}

The reduction has to approximate every vector of evaluations $Ay$. A small Euclidean residual
need not retain the moment estimate that made the first approximation possible. The following
lemma produces a residual of the form $A^+y^+$ and controls both the lengths of the new rows and
a higher moment of their evaluations.

\begin{lemma}[Interlaced projections]\label{det:interlace}
Let $A$ have rows $a_1,\ldots,a_M$ in a Euclidean space, with $|a_i|\le s$ for some $0<s\le1$.
Fix $q\ge2$ and $R\ge1$. A test is any vector $y$ for which $|y|\le1$ and $\|Ay\|_q\le R$.
Fix an integer $J\ge1$, a tolerance $0<u\le s$ and a discarded
fraction $0<\omega\le1/2$. Write
$$
 L=\max\{2,s^2/(4u)\},\qquad L_*=2^{2J-1}L,\qquad
 r=8qJ,\qquad P=40J\omega^{-2}L_*R .
$$
There exist a set $I$ of at least $(1-\omega)M$ labels, a symmetric
matrix $B$ on $I$, and rows $a_i^+$ on $I$ satisfying
$$
       |B_{ij}|\le s^2,\qquad |a_i^+|\le\sqrt3\,s/J,\qquad
       \max_{i\in I}\|B_{i,\cdot}\|_r\le P .                  \eqid{eq:residual-bounds}
$$
For each test $y$ there is a list of at most $s^2/u^2$ indices in $I$, with repetitions allowed,
coefficients $c_o$ for the individual occurrences, and a vector $y^+$ such that
$$
  |y^+|\le1,\qquad \|A^+y^+\|_{2r}\le R^+:=4P/u,
 \eqid{eq:residual-moment}
$$
$$
 \|(Ay)_I-Bc-A^+y^+\|_\infty\le 5u/2,\qquad
   \sum_o c_o^2\le s^{-2},\qquad \sum_o|c_o|\le u^{-1}.
 \eqid{eq:list}
$$
In the product $Bc$ we use $c_j=\sum_{o:j(o)=j}c_o$, where $j(o)$ is the index at occurrence
$o$. The two sums in \eqref{eq:list} are taken before repeated indices are combined. The rows of
$A^+$ lie in a common Euclidean space. The set $I$ and the matrices $B,A^+$ are the same for all
tests.
In the graph joining $i\ne j$ when $|B_{ij}|\ge b$, every closed neighborhood has at most
$1+(P/b)^r$ vertices for each $b>0$. Moreover,
$$
 R^+\le 320J\omega^{-2} 2^{2J} R u^{-2},\qquad
 \log(2R^+)\le \log(2R)+C\{J+\log(2/u)+\log(2/\omega)\}.    \eqid{eq:moment-growth}
$$
\end{lemma}

\noindent\textit{Proof of Lemma~\ref{det:interlace}.} We first prove the assertion for $s=1$.
Write $e=u$ and $L_m=2^mL$
for $0\le m<2J$.

\subsection*{Balancing the distances}

Take $Q_m$ with $1/Q_m=1/(2q)-m/(8qJ)$ for $0\le m<2J$.
Thus $2q\le Q_m<4q$. For $\theta\in[0,1]^M$ let $H_\theta$
denote the matrix with rows $\theta_j a_j$, and define
$$
 K_m(\theta)=\{x:|x|\le1,\ \|H_\theta x\|_{Q_m}\le R\}.
$$
These bodies increase with $m$, and each contains every test, since $\|H_\theta
y\|_{Q_m}\le\|Ay\|_q$. For a symmetric body $C$ write $h_C(v)=\sup_{x\in C}|\langle
v,x\rangle|$.
Put $D_0=4/\omega$ and choose $\theta$ minimizing
$$
 \frac1{2J}\sum_{m=0}^{2J-1}\sum_i
     \tfrac12\operatorname{dist}(a_i,L_m K_m(\theta))^2
        +D_0\sum_j(1-\theta_j).                                \eqid{eq:balance}
$$
Decreasing $\theta_j$ enlarges all the moment balls at once. The linear
penalty limits the number of coordinates on which this enlargement is
substantial. The distances are continuous in $\theta$. For the lower
semicontinuity, projections have norm at most $|a_i|\le1$, and limits
of feasible points remain feasible. For the upper semicontinuity, a
contraction of a feasible point by $1-\epsilon$ satisfies the moment
constraint strictly and is feasible for all nearby parameters. Letting
$\epsilon\downarrow0$ proves continuity. The minimum in \eqref{eq:balance}
is therefore attained. Comparison with $\theta=(1,\ldots,1)$ gives
$\sum_j(1-\theta_j)\le\omega M/8$. Thus at most $\omega M/4$ indices
lie outside $I_0=\{j:\theta_j\ge1/2\}$.

Write $H_0=H_\theta$, $v_{im}=P_{L_m K_m(\theta)}a_i$, and
$\delta_{im}=\langle a_i-v_{im},v_{im}\rangle$.

Projection onto the moment sublevel alone has length at most one,
because the origin belongs to that convex sublevel. Thus the ball of
radius $L_m\ge2$ is inactive. Write $J_Q(x)_j=|x_j|^{Q-2}x_j$.
There is a unique $\lambda_{im}\ge0$ such that
$$
 a_i-v_{im}=\lambda_{im}H_0^{\top}J_{Q_m}(H_0v_{im}),
 \qquad \lambda_{im}=0\quad\hbox{if }\|H_0v_{im}\|_{Q_m}<L_mR.
$$
The strict feasibility of zero gives existence. At an active point the gradient of the constraint
$\|H_0v\|_{Q_m}^{Q_m}/Q_m$ paired with $v_{im}$ equals $(L_mR)^{Q_m}>0$, giving
uniqueness. In particular
$$
        0\le\lambda_{im}\le 1/[4(L_mR)^{Q_m}],
        \qquad 0\le\delta_{im}\le1/4 .                      \eqid{eq:multiplier}
$$
Indeed $\lambda_{im}(L_mR)^{Q_m}=\delta_{im}$ on the active
constraint, and completing a square bounds the latter quantity.

The estimate for the columns comes from differentiating the objective in \eqref{eq:balance}. The
$\theta_j$-derivative of $\tfrac12\operatorname{dist}(a_i,L_mK_m(\theta))^2$ is
$$
        \lambda_{im}\theta_j^{Q_m-1}|a_j\cdot v_{im}|^{Q_m} .  \eqid{eq:derivative}
$$
For a direct verification fix $i,m$, write
$g_\theta(x)=(\|H_\theta x\|_{Q_m}^{Q_m}-(L_mR)^{Q_m})/Q_m$
and use the Lagrangian $|a_i-x|^2/2+\lambda g_\theta(x)$.
Let $V(\theta)$ be the minimum distance value and write $v,\lambda$
and $v',\lambda'$ for the solutions at two parameters $\theta,\theta'$.
The saddle inequalities give
$$
 \lambda[g_{\theta'}(v')-g_\theta(v')]
 \le V(\theta')-V(\theta)
 \le\lambda'[g_{\theta'}(v)-g_\theta(v)].
$$
The projections converge to $v$ as $\theta'\to\theta$ by compactness, continuity of the value and
uniqueness of the minimizer. Equation \eqref{eq:multiplier} bounds the
multipliers, and the normal equation and complementarity identify
each limit with $\lambda$. Difference quotients now give \eqref{eq:derivative},
including the inward one-sided derivatives at the faces of the cube.
If $\theta_j>0$, decreasing this coordinate is an admissible variation at the minimum. Formula
\eqref{eq:derivative} bounds the sum of its nonnegative contributions over all $i,m$ by $2JD_0$.
Multiplying by $\theta_j$ and keeping just one $m$ gives
$$
       \sum_i\lambda_{im}|(H_0v_{im})_j|^{Q_m}
             \le 2J D_0\theta_j\le2J D_0.                     \eqid{eq:columns}
$$
The left side vanishes when $\theta_j=0$.

\subsection*{Comparing two projected families}

Form the matrix $E_m$ of normal coefficients. Its $i$th row is $\lambda_{im}J_{Q_m}(H_0v_{im})$,
so the normal equation reads
$$
           (\langle v_{im},x\rangle)_i=Ax-E_mH_0x .            \eqid{eq:normal}
$$
Both the rows and the columns of $E_m$ have $\ell_{Q_m'}$ norm at most
$e_m=8J/(\omega L_mR)$, with $Q_m'=Q_m/(Q_m-1)$.
For rows this follows from \eqref{eq:multiplier} and $\|H_0v_{im}\|_{Q_m}\le L_mR$.
For a fixed column, the $Q_m'$th power of its norm is
$\sum_i\lambda_{im}^{Q_m'}|(H_0v_{im})_j|^{Q_m}$. Bound $\lambda_{im}^{Q_m'-1}$ by
\eqref{eq:multiplier} and apply \eqref{eq:columns} to the remaining sum. The column norm is at
most
$$
  4^{-1/Q_m}(2JD_0)^{1-1/Q_m}/(L_m R)\le e_m.
$$
Interpolating the maps $E_m:\ell_1\to\ell_{Q_m'}$ and
$E_m:\ell_{Q_m}\to\ell_\infty$ gives
$$
 \|E_m z\|_b\le e_m\|z\|_a\qquad
 (1\le a\le Q_m,\quad 1/b=1/a-1/Q_m).                        \eqid{eq:interpolation}
$$
With input exponent $q$, the output exponent satisfies $b\le2q$. The monotonicity of
counting-measure norms, $\|H_0y\|_q\le R$ and \eqref{eq:normal} therefore show that the vector
$(\langle v_{im},y\rangle)_i$ has $\ell_{2q}$ norm at most $(1+8J/\omega)R$.

For $k<m$, define $T^{km}_{ij}=\langle v_{ik},v_{jm}\rangle$ on
$I_0\times I_0$. For $x=v_{ik}$, the first term of \eqref{eq:normal} on $I_0$
has $\ell_{Q_k}$ norm at most $2L_k R$ since $\theta_j\ge1/2$
there. On the full input one has $\|H_0v_{ik}\|_{Q_k}\le L_kR$.
Apply \eqref{eq:interpolation} to $H_0v_{ik}$ with input exponent $Q_k$. The reciprocal gap
$1/Q_k-1/Q_m\ge1/(8qJ)$ gives output exponent at most $r$. Since $Q_k\le r$, every row of $T^{km}$
has $\ell_r$ norm at most
$$
       (2+e_m)L_kR\le (2+8J/\omega)L_*R=:P_0 .
$$
Only the output $Ax$ has been restricted to $I_0$. The estimate for $E_mH_0x$ is on the full set of indices.

Let $U,V$ have as their $i$th rows the averages of $v_{im}$ over the
even and the odd indices $m$, respectively. Put $\overline V=(U+V)/2$,
$D=(U-V)/2$, and denote their rows by $w_i,d_i$. The identity
$$
 \overline V\,\overline V^{\top}=B+DD^{\top},
 \qquad B=\tfrac12(UV^{\top}+VU^{\top})
$$
separates the mixed products of the two families from their squared difference.

Writing $T=(2J^2)^{-1}\sum_{e,o}T^{\min(e,o),\max(e,o)}$ over
all even--odd pairs gives $B=T+T^{\top}$. Each row of $T$ on
$I_0$ has $\ell_r$ norm at most $P_0/2$. Exclude columns whose
$\ell_r$ norm exceeds $(4/\omega)^{1/r}P_0/2$. Summing the
$r$th powers over both indices excludes at most $\omega M/4$ labels.
On the remaining set $I$ the row and column norms of $B$ are at
most $P$ and $|B_{ij}|\le1$. The two discarded fractions sum to
at most $\omega$.

\subsection*{Separation of the two means}

For brevity write $K_m=K_m(\theta)$. The projection identity and the inclusions
$L_{m-1}K_{m-1}\subset\tfrac12 L_mK_m$ give
$$
 L_m h_{K_m}(a_i-v_{im})=\delta_{im},\qquad
 \sum_m\delta_{im}\le1 .                                    \eqid{eq:projection-energy}
$$
To prove the sum estimate, put $v_{i,-1}=0$. Since $v_{i,m-1}\in\frac12 L_mK_m$ and $a_i-v_{im}$
supports $L_mK_m$ at $v_{im}$, one has
$$
 2\langle a_i-v_{im},v_{im}-v_{i,m-1}\rangle\ge\delta_{im}.
$$
Expanding $|a_i-v_{i,m-1}|^2$ around $v_{im}$ shows that the decrease in squared distance
dominates $\delta_{im}+|v_{im}-v_{i,m-1}|^2$. Summation telescopes from $|a_i|^2\le1$. A test
belongs to all the $K_m$, so
$|\langle a_i-w_i,y\rangle|\le1/(4L)\le e$.
Let $b_j=v_{i,2j+1}-v_{i,2j}$. Nested projection gives
$\sum_j|b_j|^2\le1$. For $j<k$, nested bodies and \eqref{eq:projection-energy} give
$$
 |\langle b_j,b_k\rangle|
 \le (L_{2j}+L_{2j+1})\left(\delta_{i,2k}/L_{2k}+
                                  \delta_{i,2k+1}/L_{2k+1}\right).
$$
Summing first in $j<k$ gives a factor at most $L_{2k}$,
then summing in $k$ gives one. Consequently
$$
 |d_i|\le\sqrt3/(2J),\qquad
 |w_i|^2+|d_i|^2=(|U_i|^2+|V_i|^2)/2\le1.                 \eqid{eq:separation}
$$

\subsection*{An exact constrained pursuit}

Write $G_i=(w_i,d_i)$ for $i\in I$. Then $GG^{\top}=B+2DD^{\top}$
and $|G_i|\le1$. Put $T_0=P/e$ and
$C=\{x:\|Gx\|_{2r}\le T_0\}$. Starting at $x_0=(y,0)$ gives
$\|Gx_0\|_{2r}\le(1+8J/\omega)R\le T_0$ by \eqref{eq:normal}--\eqref{eq:interpolation}.
At any label $j$ with $|G_j\cdot x|>e$ record
$c=G_j\cdot x$, take $x'=x-cG_j$ and put $x_+=P_Cx'$.
The first operation decreases $|x|^2$ by $(2-|G_j|^2)c^2\ge c^2$. For the projection, zero is
strictly feasible, so there is a vector $b'$ satisfying
$$
 x'-x_+=G^{\top}b',\qquad
 b'\cdot Gx_+=T_0\|b'\|_{(2r)'} .
$$
At a boundary point $b'$ is a nonnegative multiple of $J_{2r}(Gx_+)$, while at an interior point
one takes $b'=0$. The second operation consequently satisfies $|x'|^2-|x_+|^2=2\langle
b',Gx_+\rangle+|x'-x_+|^2\ge2T_0\|b'\|_{(2r)'}$. Thus the process terminates after at most
$e^{-2}$ occurrences. The sum $\beta$ of all the normal coefficients satisfies
$\|\beta\|_{(2r)'}\le (2T_0)^{-1}$ and the coefficients obey \eqref{eq:list}
by the decrease in $|x|^2$ and Cauchy--Schwarz. The terminal $x$ satisfies $|x|\le1$ and
$\|Gx\|_\infty\le e$. Its second Euclidean component is $-y^+$, where
$y^+=D^{\top}(c+\beta)$ and $|y^+|\le1$.
At termination $x=(y,0)-G^{\top}(c+\beta)$, with $c$ the aggregated list. Multiplying by $G$ and
using $GG^{\top}=B+2DD^{\top}$ gives
$$
                  \overline V y=Bc+2Dy^++Gx+B\beta .
$$
By symmetry and the row bound, the norms of $B:\ell_1\to\ell_r$ and
$B:\ell_{r'}\to\ell_\infty$ are at most $P$. Interpolation yields
$\|B\|_{(2r)'\to 2r}\le P$. Also
$\|B\beta\|_\infty\le P\|\beta\|_{r'}\le P/(2T_0)\le e/2$.
Put $a_i^+=2d_i$. Equations \eqref{eq:projection-energy}--\eqref{eq:separation} and the terminal identity
prove the row-radius bound and \eqref{eq:list}. In the same identity
$\|\overline V y\|_{2r}\le T_0$, $\|Bc\|_{2r}\le P/e=T_0$,
$\|Gx\|_{2r}\le T_0$, and $\|B\beta\|_{2r}\le P/(2T_0)\le T_0$.
This gives \eqref{eq:residual-moment}.

For $s<1$, apply the unit-radius case to $\widehat A=A/s$ with
$\widehat R=R/s$ and $e=u/s$. Its constant $\widehat P$ satisfies
$s^2\widehat P\le P$, since
$s\max\{2,s/(4u)\}\le\max\{2,s^2/(4u)\}$.
Use $B=s^2\widehat B$, $c=\widehat c/s$ and
$A^+=s\widehat A^+$. Multiplying the scalar identity by $s$ gives
\eqref{eq:list}, and the moment of the new term is at most
$4s\widehat P/e=4s^2\widehat P/u\le4P/u$. This proves \eqref{eq:residual-bounds}--\eqref{eq:moment-growth}
at every radius.\hfill\qedsymbol

\section{Counting localized approximations}
\label{sec:support-count}

At each use of Lemma~\ref{det:interlace}, a test contributes a short sum of columns. Their large
entries lie in graph neighborhoods. The graphs at later stages may have very large degrees, but
the lists of columns at those stages are shorter. A count using the largest degree throughout
would lose this compensation. We retain only the part of each partial sum that exceeds a
threshold before adjoining the preceding layer.

\begin{lemma}[Localized approximations on successive graphs]\label{det:count}
Fix $m\ge1$, a set of $M$ labels and graphs $G_1,\ldots,G_m$ on
these labels. Their maximum closed-neighborhood sizes are bounded
by $1\le D_1\le\cdots\le D_m$. Fix nonnegative integer length
bounds $K_j$, let $K=\max(1,\sum_jK_j)$, and let each occurrence
in layer $j$ carry a coefficient from a fixed alphabet of size at
most $Q$, where $Q\ge1$. At each layer join two occurrences whose
labels coincide or have graph distance at most two in $G_j$.
For a resulting component $L$ write
$$
 U_L=\bigcup_{a\in L}N_{G_j}(a).
$$
Its labels and coefficients determine a scalar function $P_L$
supported in $U_L$. The $U_L$ at one layer are pairwise disjoint.
Fix $0<h<1$ and put $\epsilon=2^{-m-1}h$. An allowed test $z\in\mathbb R^M$ satisfies $\sum_i|z_i|^p\le W$,
where $p>0$, and admits lists of lengths at most $K_j$ with
$$
 \left\|z-\sum_{j=1}^m\sum_{L\ {\rm in\ layer}\ j}P_L\right\|_\infty
       \le\epsilon/2.                                             \eqid{eq:lists-error}
$$
For instance, $P_L$ may be the rounded sum of matrix columns on
$U_L$, provided the rounding, the omitted entries and the analytic
approximation together yield \eqref{eq:lists-error}.

Define
$$
 L_j^*=\log\!\left[
 \frac{64m}{h}(2+K)(2+Q)(2+D_j)
             \big(2+W(2/\epsilon)^p\big)\right].
$$
There is a collection of scalar approximants indexed by nonzero
profiles $\nu=(\nu_1,\ldots,\nu_m)$, where $\nu_j$ counts the
occurrences from layer $j$. For each profile its size is at most
$$
 M\exp\!\left(C_m\sum_j\nu_j L_j^*\right).                      \eqid{eq:count}
$$
Every allowed $z$ admits members of the collection with disjoint
supports whose sum approximates $z$ in $\ell_\infty$ within $h$.
For these members the sum of the profiles is coordinatewise at
most $(K_1,\ldots,K_m)$.

\end{lemma}

\begin{proof}
To construct the approximations, proceed from $j=m$ down to $1$ and put
$\theta_j=2^{m-j}\epsilon$. At the start of layer $j$ the surviving
records from the deeper layers have disjoint supports $E_H$. Each
record determines both $E_H$ and its scalar approximant there. Join
each component $L$ of layer $j$ to every deeper record $H$ for
which $U_L\cap E_H\ne\varnothing$. Include isolated vertices.
Within every connected component of this bipartite graph add all
approximants $P_L$ and all deeper approximants, extending each by zero.
The resulting component sums have disjoint supports. Retain each
sum only where its modulus is at least $\theta_j$, and discard it
if its retained support is empty. Its new record consists of the
constituent records and the layer-$j$ occurrences. Thus its
retained support and its values can be recovered from its record.

If $F_j=\sum_L P_L$ is the sum of the original layer-$j$ approximants
and $T_j$ is the sum immediately after clipping layer $j$, the
construction gives $T_{m+1}=0$ and
$$
 T_j=(F_j+T_{j+1})\boldsymbol1_{\{|F_j+T_{j+1}|\ge\theta_j\}},
 \qquad
 \left\|T_j-\sum_{r\ge j}F_r\right\|_\infty
       \le\sum_{r\ge j}\theta_r .
$$
In particular \eqref{eq:lists-error} gives $\|T_1-z\|_\infty<h$. At a retained index outside the
union of all neighborhoods from layers $r<j$,
\eqref{eq:lists-error} and the losses at deeper layers give
$$
 |z_i|\ge\theta_j-\sum_{r>j}\theta_r-\epsilon/2=\epsilon/2 .
$$
Consequently every actual retained support at layer $j$ has size
at most
$$
 R_j=\left\lceil W(2/\epsilon)^p\right\rceil
                         +\sum_{r<j}K_rD_r.                       \eqid{eq:retained}
$$
We may discard a record as soon as its support is empty or exceeds the bound
\eqref{eq:retained}. This decision depends only on the record. The preceding estimate shows that
every nonempty component arising from an admissible test survives.

Fix a label $i$. We claim that the number of layer-$j$ records with profile $\nu$ whose retained
support contains $i$ is at most
$$
 \exp\!\left(A_{m-j+1}\sum_{r\ge j}\nu_r L_r^*\right),          \eqid{eq:root}
$$
where the constants depend only on the indicated depth. A single
layer-$j$ component $L$ of $u$ occurrences for which an address
$i\in U_L$ is given has at most
$$
        \big[(2+K)^3(2+Q)(2+D_j)^2\big]^u
$$
encodings. Choose a first anchor from $N_{G_j}(i)$, a rooted
spanning tree of the occurrences, and each subsequent anchor
within graph distance two of its parent. Include its rounded
coefficient at each vertex. The labels need not be distinct.
This proves \eqref{eq:root} for the deepest layer.

For $j<m$, root a spanning tree of the bipartite component at a vertex
whose set contains $i$ and order its children. If the root is a
layer-$j$ component $L$, its choices are bounded by the preceding
single-layer count. If it is a deeper record $H$, apply
\eqref{eq:root} at layer $j+1$. Each edge now requires one choice
of an index in the intersection of the sets at its two ends.

Suppose first that the parent is a deeper record $H$ and its child
is a component $L$ at layer $j$. Once $H$ is known, the intersection
index can be chosen from its retained support, in at most
$$
 R_{j+1}\le 1+W(2/\epsilon)^p+KD_j
$$
ways. With this index fixed, the single-layer count applies to $L$.
The logarithm of the extra factor is bounded by a constant times
$L_j^*$, independently of all deeper degrees. In the other direction,
when $L$ is the parent, at most $KD_j$ indices lie in $U_L$.
For each such choice, \eqref{eq:root} bounds the number of possible
children $H$. Each contains an occurrence at a layer $r>j$, and
$D_j\le D_r$, so its exponent absorbs the factor $KD_j$.

It remains to count the tree and the distribution of the occurrences
among its vertices. If $u=\sum_{r\ge j}\nu_r\le K$, number the $u$
occurrences and assign to each a layer and one of at most $u$
vertices. Specify the type of each vertex, the root, its parent
relations and the order of the children. Together these choices
contribute at most $\exp(C_m u\log(2+K))$. The exponents in the
counts for the children add over disjoint occurrence sets. Each
nonroot vertex incurs the cost of just one intersection index, so
induction proves \eqref{eq:root}, with $A_d=A_{d-1}+C_m$ after
enlarging $C_m$. Finally, a nonempty support contains some $i$.
Summing the bound for fixed $i$ over the $M$ possible initial indices
gives \eqref{eq:count}.

\end{proof}

\subsection*{The constraints at a fixed height}

For $p>2$ and a dyadic height $0<t\le1$, assume $h\le t/16$. Put
$f_t(0)=0$ and
$$
 f_t(u)=|u|^p\big(1-|\log_2(|u|/t)|\big)_+\quad(u\ne0),\qquad
 F_t(z)=\sum_{t/4\le|z_i|\le3t}|z_i|^p.
$$
Apply $f_t$ to each final approximant on its support and impose the
resulting scalar constraints with squared thresholds
$\lambda_\nu^2=C'_m\sum_j\nu_j L_j^*$. With $C'_m$ sufficiently
large, the sum of $\exp(-\lambda_\nu^2/16)$ over these constraints
is at most $Mh^2/256$. Indeed, after \eqref{eq:count}, the contribution for
each nonzero profile is at most $M\exp(-2\sum_j\nu_j L_j^*)$.
For total size $u\ge1$ there are at most $m^u$ profiles, and
$L_j^*\ge\log(64m/h)$. The resulting bounds are summable over dyadic heights when $h\le t/16$.

For the disjoint approximants of an actual test write $v_C$ for
their band vectors. Each nonzero entry occurs at a label with
$t/4\le |z_i|\le3t$ and has modulus at most $(2t)^p$. Consequently
$$
 \sum_C\|v_C\|_2^2\le C_p t^p F_t(z),\qquad
 \left\|f_t(z)-\sum_Cv_C\right\|_1\le C_p(h/t) F_t(z).
$$
The second estimate uses the $C_p t^{p-1}$ Lipschitz bound for
$f_t$. For any coloring $\|\chi\|_\infty\le1$ satisfying the
scalar constraints, Cauchy--Schwarz and the profile budget yield
$$
 |\langle\chi,f_t(z)\rangle|
 \le C_p(h/t)F_t(z)+
 C_{p,m}t^{p/2}\sqrt{F_t(z)}
                  \left(\sum_j K_j L_j^*\right)^{1/2}.           \eqid{eq:band-bound}
$$

\section{Discrepancy across all heights and resolutions}\label{full:halving}

Fix $p>2$ and $\sigma>0$. Choose
$$
 0<a<\min\{1/4,(p-2)/2,\sigma/2\},\qquad d_*=\lceil 1/a\rceil .
$$
All constants in this section may depend on $p,\sigma$ through these
choices. Use the normalized rows, the number $W\ge1$ and the body
$K$ of Section~\ref{full:normalization}. The scalar band at height
$t=2^{-\ell}$, $\ell\ge0$, is the function $f_t$ defined in
Section~\ref{sec:support-count}. We use the relative
resolutions $\tau_v=2^{-v}/256$ for $0\le v\le v_*$, with $v_*$
chosen below. At a given height and resolution write
$$
 h=\tau_vt,\qquad u=c_*h,\qquad J=\lceil 2u^{-a}\rceil ,
 \qquad s_j=u^{aj}\quad(0\le j\le d_*) .
$$
Take $c_*>0$ so small that $(5d_*+1)c_*\le2^{-d_*-2}$.
We construct $d_*$ successive approximations at this pair $(v,\ell)$.

Start at $s_0=1$, $q_0=p$, $R_0=W^{1/p}$ with the original rows on
all $M$ labels. At stage $j$ use the projection and pursuit construction
of Section~\ref{sec:weighted-interlacing} at absolute tolerance $u$,
with row norm bound $s_j$ and the discarded fraction
$$
 \omega=\frac{c}{d_*(1+\ell)^2(1+v)^2}.
$$
Here $c$ is absolute and sufficiently small. Since $j<d_*$, $s_j>u$.
The resulting rows for the next stage have length at most
$\sqrt3\,s_j/J\le s_{j+1}$. They are defined on the locally retained
labels and are independent of the test. On this same set the new test
has norm at most one and its evaluations have $\ell_{q_{j+1}}$ norm
at most $R_{j+1}$, with $q_{j+1}=16Jq_j$. Thus the construction can
be repeated with all unit tests satisfying that moment bound.
After the last stage the remaining evaluations are bounded by
$s_{d_*}\le u$ at every retained label.

Let $B_j$ be the symmetric matrix at stage $j$ and $P_j$ the bound on its $\ell_{8q_jJ}$ row
norms. Join distinct retained indices $i,k$ when $|(B_j)_{ik}|\ge u^2$,
and take closed neighborhoods. Unretained labels may be included in
the graph as isolated vertices. The bounds \eqref{eq:residual-bounds}--\eqref{eq:moment-growth}, with ratio
two between successive dilates, imply
$$
 \begin{aligned}
  q_j&=(16J)^jp,\\
  \log(2+P_j)+\log(2+R_{j+1})
     &\le C_{p,\sigma}\{\log(2+W)+J+\log(2/u)+\log(2/\omega)\}.
 \end{aligned}
$$
Indeed each step increases $\log(2+R_j)$ by at most an absolute
multiple of $J+\log(2/u)+\log(2/\omega)$, and the depth is fixed.
There are at most $1+(P_j/u^2)^{8q_jJ}$ labels in any closed
neighborhood. Since $\ell+v\le C\log(2/u)$ and
$\log(2/u)\le C_a u^{-a}$, the closed degrees have nondecreasing
upper bounds $D_j$ satisfying
$$
 \log(2+D_j)\le
  C_{p,\sigma} u^{-a(j+1)}\{\log(2+W)+u^{-a}\}.                \eqid{eq:degrees}
$$
The number of recorded coefficients at stage $j$ is at most
$s_j^2/u^2$. Each has modulus at most $s_j^{-1}\le u^{-1}$
by \eqref{eq:list}. Round each towards zero on the mesh $u^3$, and retain each
matrix column only on the closed neighborhood of its label. At each
output coordinate, rounding costs at most
$u^3(s_j^2/u^2)s_j^2\le u$. Omitting the small entries costs at
most $u^2\sum_o|c_o|\le u$. By \eqref{eq:list}, summing these rounded, truncated
columns over all stages approximates the original evaluations with
uniform error at most $(5d_*+1)u$ on the final locally retained set.
The coefficient alphabet has at most $1+2u^{-4}$ elements.

Let $I$ be the intersection of the final retained sets over all heights
and resolutions. The number of excluded labels is at most
$M\sum_{v,\ell}d_*\omega\le M/128$. In every collection use only $I$ as
output labels. The earlier stages can use any anchor that was locally
retained at that stage. Its neighborhood is known on the original set
of $M$ labels and obeys \eqref{eq:degrees}. Equivalently, extend the evaluation vector
and all predicted columns by zero off $I$. These conventions preserve
the estimates for the normal matrices, for which all input labels at
the given stage were used.

\subsection*{Two consecutive resolutions}

For each fixed test, height, and resolution choose one of its lists
at every stage, using the same lists whenever that resolution recurs.
Write $z_i=a_i(y)$ on $I$ and $z_i=0$ elsewhere. At resolution $v$,
form the successive clipped sums $T_j$ as in
Section~\ref{sec:support-count}, using $m=d_*$ stages with
its local $h=\tau_vt$. The choice of $c_*$ gives \eqref{eq:lists-error} with
$\epsilon=2^{-d_*-1}h$. Denote the final scalar approximant by $Z_v$,
so that $\|Z_v-z\|_\infty\le\tau_vt$.

For $v\ge1$ we estimate $f_t(Z_v)-f_t(Z_{v-1})$. The summability in $v$ depends on counting the
two approximations together, even though their graphs and chosen indices can differ.
At stage $j$ take the union of the two graphs and mark each occurrence with its resolution. Its
column keeps its own original truncation.
Carry the two scalar predictions as an ordered pair. At each clipping
step apply the respective scalar threshold to each coordinate, and
retain the union of the two resulting supports. Since clipping acts
pointwise, each coordinate of the resulting sum is exactly the
single-resolution approximant already defined. At an output label of
a retained record outside all earlier neighborhoods, the proof of
\eqref{eq:retained} for its nonzero coordinate gives $|z_i|\ge 2^{-d_*-2}\tau_vt$.
Thus \eqref{eq:retained} holds with the smaller resolution $h=\tau_vt$ and the
combined lists and graphs. The rooted count \eqref{eq:count}--\eqref{eq:root} is unchanged
for these pairs, since each record determines its support and both approximating functions. For
$v=0$ use the single collection and the difference from
zero, so $f_t(Z_{-1})=0$.

On a final record $C$ use as the scalar constraint the difference of the
two band approximants on its retained support, extended by zero elsewhere. These
vectors $b_C$ have disjoint supports for any actual pair of lists.
Put $F_t(y)=\sum_{i\in I:\ t/4\le |z_i|\le3t}|z_i|^p$. Since $f_t$
is $C_pt^{p-1}$-Lipschitz and supported where $t/2\le|z|\le2t$,
$$
       \sum_C |b_C|^2\le C_p\tau_v^2t^pF_t(y).                 \eqid{eq:paired-band}
$$
Indeed, for $v\ge1$ the pointwise difference is at most $3C_p\tau_vt^p$ because
$\tau_{v-1}=2\tau_v$. The number of its possible nonzero coordinates is at most $(4/t)^pF_t(y)$,
since each such coordinate has
$t/4\le|z_i|\le3t$. For $v=0$ the same bound holds after enlarging
the absolute factor since $\tau_0=1/256$.

Use the thresholds leading to \eqref{eq:band-bound}. For the paired lists let $K_j$ be the sum
of the two length bounds at stage $j$. For the pair, if $u=c_*\tau_vt$, then
$$
 K_j\le2u^{-2+2aj},\qquad Q\le 2(1+2u^{-4}),\qquad 0\le j<d_* .
$$
The degrees of the union obey \eqref{eq:degrees} with a larger constant, because the
two tolerances differ by a factor two. In the notation $L_j^*$ of the
collection, $\log(2+\sum K_j)$, $\log(2+Q)$, and
$\log(2+W(2/\epsilon)^p)$ together are bounded by
$C_{p,\sigma}[\log(2+W)+\log(2/u)]$. Hence
$$
 \sum_{j=0}^{d_*-1}K_j L_j^*
 \le C_{p,\sigma}u^{-2-2a}\log(2+W).                           \eqid{eq:entropy}
$$
For the degree term, multiplication in \eqref{eq:degrees} gives the powers
$u^{-2+aj-a}$ and $u^{-2+aj-2a}$, whose largest values occur at
$j=0$. This estimate is where the decrease in the lengths of the
successive lists enters.

For fixed $v$ and $t$, \eqref{eq:count} bounds the sum of $\exp(-\lambda^2/16)$ by
$M(\tau_vt)^2/256$ after the constant in the squared thresholds is chosen. Its sum over all
$v,\ell$ is less than $M/64$. Impose $\chi_i=0$ for each $i\notin I$ and impose orthogonality to
the $s$ cost vectors, all with threshold zero. Since $|I^c|,s\le M/128$, the combined sum is
below $M/16$. Lemma~\ref{det:LM} now gives one coloring $\chi$ with at least $M/2$ saturated
coordinates. For any test, disjointness, \eqref{eq:paired-band},
Cauchy--Schwarz and the profile budget \eqref{eq:entropy} show that
$$
 |\langle\chi,f_t(Z_v)-f_t(Z_{v-1})\rangle|
 \le C_{p,\sigma}\tau_v^{-a}t^{p/2-1-a}
                     \sqrt{F_t(y)\log(2+W)}.                   \eqid{eq:paired-discrepancy}
$$
The total $\ell_1$ error at the terminal resolution is bounded by
$C_p\tau_{v_*}F_t(y)$, as in \eqref{eq:band-bound}.

For $|z_i|\le1$ one has $\sum_{\ell\ge0}f_{2^{-\ell}}(z_i)=|z_i|^p$,
and $\sum_tF_t(y)\le C F(y)$ by bounded overlap. The exponent $\kappa=p/2-1-a$ is positive.
Cauchy--Schwarz yields $\sum_t t^\kappa\sqrt{F_t(y)}\le(\sum_t
t^{2\kappa})^{1/2}(\sum_tF_t(y))^{1/2}\le C_{p,a}\sqrt{F(y)}$, while $\sum_{v=0}^{v_*}\tau_v^{-a}
\le C_a\tau_{v_*}^{-a}$. Summing \eqref{eq:paired-discrepancy} and the terminal errors gives
$$
 \left|\sum_i\chi_i|a_i(y)|^p\right|
 \le C_{p,\sigma}\left[\tau_{v_*}F(y)+
       \tau_{v_*}^{-a}\sqrt{F(y)\log(2+W)}\right],\quad y\in K.
$$
Choose $v_*=\lceil\log_2(W)/(2+\sigma)\rceil$. On $\partial K$,
divide by $F(y)=W$ and use $\tau_{v_*}\asymp W^{-1/(2+\sigma)}$.
Since $2a<\sigma$, the resulting one-step relative error is at most
$$
                 C_{p,\sigma}W^{-1/(2+\sigma)}.                 \eqid{eq:step}
$$
In detail, $W^{a/(2+\sigma)-1/2}\sqrt{\log(2+W)}
\le C_{a,\sigma}W^{-1/(2+\sigma)}$ for $W\ge1$.
Homogeneity extends \eqref{eq:step} from $\partial K$ to every test.

\section{Iteration and geometric consequences}\label{full:conclusion}

\begin{proof}[Completion of Proposition~\ref{det:reduction}]
Choose the sign of $\chi$ that gives at least $M/4$ negative
saturated coordinates. The nonnegative multipliers $1+\chi_i$
preserve the total cost. Absorbing each surviving multiplier into the
cost and its $p$th root into the corresponding form reduces the
number of labels to at most $3M/4\le27N/32$. Undo the common
normalization of the forms and add the weights of its copies with the factors $1/k_i$ from the
splitting. Each original label with positive weight has at least one surviving copy. Estimate
\eqref{eq:step} therefore proves the proposition.
\end{proof}

\begin{proof}[Proof of Theorem~\ref{full:rows}]
Repeat the reduction with the new weighted family, stopping at
$$
 N\le K_*=\max\{B_{p,\sigma}\eta^{-(2+\sigma)}n^{p/2},128(s+1)\}.
$$
If the initial family already meets this bound, take all weights equal to one. At every reduction of positive rank,
$s\le M/128$ and $W=N/(8r^{p/2})\ge
(B_{p,\sigma}/8)\eta^{-(2+\sigma)}N/K_*$, with $r\le n$.
Put $\gamma=1/(2+\sigma)$. If $N_t$ is the last size at which a
reduction is performed, then $N_t>K_*$ and
$N_{t-j}\ge (32/27)^jN_t$. By \eqref{eq:step}, the sum of the relative
errors is bounded by
$$
 C_{p,\sigma}B_{p,\sigma}^{-\gamma}\eta
          \sum_{j\ge0}(27/32)^{j\gamma}.
$$
Choose $B_{p,\sigma}$ so that this is at most $\eta/2\le\log(1+\eta)$.
The product of the upper factors is then at most $1+\eta$.
The product of the lower factors is at least one minus the sum of
the errors and hence at least $1-\eta$.

If a
vanishing-rank step occurs, its conic replacement reproduces the zero moment exactly and ends the
iteration. Composition of multipliers gives weights on the original family.

For complex forms apply the real theorem in dimension $2n$ to
$y\mapsto\operatorname{Re}a_i(y)$ with the given real costs. Test
the real inequality at $e^{i\theta}y$ and integrate in $\theta$.
Each term averages to $|a_i(y)|^p$ times the same positive mean
of $|\cos\theta|^p$. This proves Theorem~\ref{full:rows} over both fields.
\end{proof}

\begin{corollary}\label{full:embedding}
For $p>2$, $\sigma>0$, and $0<\eta<1$, every $n$-dimensional
subspace of real or complex $L_p(\mu)$ embeds into $\ell_p^m$ with
distortion at most $1+\eta$, where
$m=\lceil C_{p,\sigma}\eta^{-(2+\sigma)}n^{p/2}\rceil$.
\end{corollary}

\begin{proof}
Fix simultaneously chosen measurable representatives of a basis and let $a(x)$ be the
vector of its values. The integral of $|a(x)|_2^p$ is finite.
Push $|a(x)|_2^p\,d\mu(x)$ forward to the Euclidean unit sphere by
$a(x)/|a(x)|_2$, omitting the zero set. Its $p$th moments give the
$p$th power of the norm in basis coordinates. Finite positive
measures approximate these moments uniformly on the Euclidean unit
sphere of coefficients, by partitioning the sphere of rows into
sets of small diameter. The norm has a positive minimum on the
coefficient sphere, so this approximation can be taken with any
prescribed relative error. Apply Theorem~\ref{full:rows} to the
finite measure, absorbing its masses into the rows. Taking each
relative error at most $\eta/16$ bounds the distortion by
$[(1+\eta/16)/(1-\eta/16)]^{2/p}\le1+\eta$.
\end{proof}

In particular, if a matrix $A$ has rank $n$, its row space has dimension $n$ and the theorem
applies on the quotient by $\ker A$. Multiplying the selected rows by $w_i^{1/p}$ gives a matrix
$\widetilde A$ for which $(1-\eta)\|Ax\|_p^p\le\|\widetilde Ax\|_p^p\le(1+\eta)\|Ax\|_p^p$ for
every $x$. The construction also gives a probability law on such weighted selections with
prescribed expectations.

\begin{corollary}[Preservation in expectation]\label{cor:unbiased}
For the family and costs in Theorem~\ref{full:rows}, there is a finitely supported probability law
on weights with $\mathbb E w_i=1$ for each $i$. Every realization satisfies both moment
inequalities, preserves the $s$ costs exactly, and has at most
$C_{p,\sigma}\eta^{-2-\sigma}n^{p/2}+128(s+3)$ nonzero entries.
\end{corollary}
\begin{proof}
Let $\mathcal C$ be the set of all weights satisfying the indicated support, moment and cost
constraints together with $\sum_iw_i=N$. It is compact, since it is closed in $[0,N]^N$. Given
$h\in\mathbb R^N$, apply Theorem~\ref{full:rows} with costs $(c_i,1,h_i)$. It gives a point of
$\mathcal C$ with $\langle h,w\rangle=\langle h,\boldsymbol1\rangle$. Strict separation places
$\boldsymbol1$ in the compact convex hull of $\mathcal C$. A finite convex combination
representing this point gives the required law.
\end{proof}

\section{Spherical moments and Euclidean sections at a fixed exponent}\label{sec:spherical}

All random vectors and subspaces below are real. Let $\sigma_d$ be uniform probability measure
on $\Sph^{d-1}=\{x\in\R^d:|x|=1\}$ and $\gamma_d$ standard Gaussian measure on $\R^d$. We write
$o_{\PP}(1)$ for a quantity tending to zero in probability. In a uniform estimate this
convention applies to the supremum of the absolute error over the indicated set.

Dvoretzky's theorem \cite{Dvo} asserts that, given $k\in\mathbb N$ and $D>1$, every normed space
of sufficiently large dimension contains a $k$-dimensional subspace at Banach--Mazur distance at
most $D$ from $\ell_2^k$. Milman's probabilistic proof \cite{Mil} led to the study of the
dimension at which a typical subspace is almost Euclidean. For classical sequence spaces, Gordon,
Gu\'edon, Meyer and Pajor \cite{GGMP} obtained estimates for the Banach--Mazur distance of random
sections to Euclidean space. Paouris, Valettas and Zinn \cite{PVZ} refined the dependence on the
prescribed distortion. The case of the cube was studied by Schechtman \cite{Sch} and Tikhomirov
\cites{Tik14,Tik18}.

For $1\le p\le\infty$ let $\ell_p^N=(\R^N,\|\cdot\|_p)$, with the usual sequence norm (the maximum
norm when $p=\infty$). Given a nonzero subspace $E\subset\R^N$, put
$$
 \Delta_p(E)=\frac{\sup\{\|z\|_p:z\in E,\ \|z\|_2=1\}}
 {\inf\{\|z\|_p:z\in E,\ \|z\|_2=1\}}.
 \eqid{dist}
$$
Thus $\Delta_p(E)\le D$ if and only if there is a scalar $a>0$ such that
$$
 a\|z\|_2\le\|z\|_p\le Da\|z\|_2\qquad(z\in E).
$$
The norm $\|\cdot\|_2$ in this definition is inherited from the ambient space. For $k=\dim E$ put
$$
 \mathcal D_p(E)=d_{\rm BM}((E,\|\cdot\|_p),\ell_2^k)
 :=\inf_T\|T\|\|T^{-1}\|,
$$
where the infimum runs over the linear isomorphisms $T:(E,\|\cdot\|_p)\to\ell_2^k$. In particular,
$1\le\mathcal D_p(E)\le\Delta_p(E)$. The distinction between geometric and isomorphic estimates
also occurs in the randomized isomorphic Dvoretzky theorem of Litvak, Mankiewicz and
Tomczak-Jaegermann \cite{LMT}. Let $\mathcal G_{N,k}$ be the Grassmannian of $k$-dimensional
subspaces of $\R^N$, equipped with its rotation-invariant probability measure. For $D>1$ and
$0<\delta<1$, let
$$
 d(N,p,D,\delta)=\max\{1\le k\le N:\PP(\Delta_p(E)\le D)\ge1-\delta,\ E\sim\mathcal G_{N,k}\}.
$$
We write $d_{\rm iso}(N,p,D,\delta)$ for the same maximum with $\Delta_p(E)$ replaced by $\mathcal
D_p(E)$. A one-dimensional subspace has both distances equal to one, so these maxima exist. If
$E\subset E'$, the supremum in \eqref{dist} increases and the infimum decreases upon passing from
$E$ to $E'$. For the Banach--Mazur distance, any Euclidean norm $q$ on $E'$ satisfying
$q\le\|\cdot\|_p\le Dq$ restricts to a Euclidean norm with the same inequalities on $E$. Hence
both distances are nondecreasing under inclusion. Taking the first $k$ columns of a
Haar-distributed orthogonal matrix produces a nested family with the required marginal
distributions. The probabilities defining the two maxima are therefore nonincreasing in $k$.

For fixed $p>2$, the critical dimension is of order $N^{2/p}$. The limiting distortion on this
scale is determined by the extrema of the empirical $p$th moment of independent uniform directions
on the sphere. For $q>0$ write
$$
 \mu_q=\frac1{\sqrt{2\pi}}\int_{\R}|t|^q e^{-t^2/2}\,dt.
 \eqid{mu}
$$
For $p>2$, if $U_1,\ldots,U_m$ are independent and uniform on $\Sph^{n-1}$, define
$F:\Sph^{n-1}\to\R$ by $F(x)=\sum_{i=1}^m|\langle U_i,x\rangle|^p$. For $0<r\le1$ write
$$
 j(r)=\lfloor1/r\rfloor,\qquad H_p(r)=j(r)r^{p/2}+(1-j(r)r)^{p/2}.
 \eqid{profile}
$$
\begin{mainresult}[Spherical localization]\label{th:fixed}
Fix $p>2$ and suppose $m/n^{p/2}\to\alpha\in(0,\infty)$. For every fixed $0<r\le1$, as $n\to\infty$,
$$
 \inf_{x\in\Sph^{n-1}} F(x)\xrightarrow{\PP}\alpha\mu_p,
 \eqid{fixmin}
$$
$$
 \sup_{\substack{x\in\Sph^{n-1}\\\max_i|\langle U_i,x\rangle|^2\le r}}F(x)
 \xrightarrow{\PP}\alpha\mu_p+H_p(r).
 \eqid{fixmax}
$$
The admissible set in \eqref{fixmax} is nonempty with probability tending to one. We assign the
restricted supremum the value zero on the event that it is empty. If $X_n\in\Sph^{n-1}$ are random
directions satisfying $\sup_{x\in\Sph^{n-1}}F(x)-F(X_n)\xrightarrow{\PP}0$, then $\max_{1\le i\le
m}|\langle U_i,X_n\rangle|\xrightarrow{\PP}1$.
\end{mainresult}
In particular, when $r=1/j$ with $j$ a positive integer, the extra term in \eqref{fixmax} is
$j^{1-p/2}$. At $r=1$ the restriction disappears and the excess is one. Comparing this maximum
with the minimum yields the inherited distortion of the random section. To optimize its
Euclidean structure we shall also need the average in Lemma~\ref{ellipsoid}.

Choose $1/p<a<1/2$ and put $t=n^{-a}$. Below this height the truncated moment concentrates
uniformly. Above it, the squared inner products have total mass at most $1+o_{\PP}(1)$. The
constants in the next two lemmas may depend on $p$ and $\alpha$.
\begin{lemma}\label{trunc}
Under the assumptions of Theorem~\ref{th:fixed},
$$
 \sup_{x\in\Sph^{n-1}}\left|\sum_{i=1}^m\min\{|\langle U_i,x\rangle|^p,t^p\}-\alpha\mu_p\right|\xrightarrow{\PP}0.
 \eqid{trunclim}
$$
\end{lemma}
\begin{proof}
Let $\sigma_n$ denote the normalized surface measure on $\Sph^{n-1}$. For $q>0$ and
$x\in\Sph^{n-1}$, rotational invariance and Gaussian polar coordinates give
$$
 c_{n,q}:=\int_{\Sph^{n-1}}|\langle u,x\rangle|^q\,d\sigma_n(u)
 =\frac{\Gamma(n/2)\Gamma((q+1)/2)}{\sqrt\pi\,\Gamma((n+q)/2)}\sim \mu_q n^{-q/2}.
 \eqid{mom}
$$
Indeed a standard Gaussian vector $g$ decomposes as $g=\|g\|_2U$ with $U$ independent of $\|g\|_2$
and uniform on the sphere. The $q$th moment of $\|g\|_2$ is $2^{q/2}\Gamma((n+q)/2)/\Gamma(n/2)$
and the $q$th absolute moment of its first coordinate is $\mu_q$ as in \eqref{mu}. Their ratio
gives \eqref{mom}.

Fix $x\in\Sph^{n-1}$ and set $Z_i(x)=\min\{|\langle U_i,x\rangle|^p,t^p\}$.
Writing $U_i=g/|g|$ and splitting according to $|g|\ge\sqrt{n/2}$ gives,
for $0<s<1$,
$$
 \PP\{|\langle U_i,x\rangle|>s\}\le 2e^{-ns^2/4}+e^{-c_1n}.
$$
Since $m$ grows polynomially and $nt^2=n^{1-2a}$, truncation changes
$mc_{n,p}$ by $o(1)$. Thus $m\E Z_i(x)\to\alpha\mu_p$, independently of
$x$. Also $\sum_i\operatorname{Var}Z_i(x)\le mc_{n,2p}\le Cn^{-p/2}$
and $|Z_i(x)-\E Z_i(x)|\le t^p$. For each fixed $\xi>0$, Bernstein's
inequality on an $n^{-p-2}$-net $\mathcal N_n$ gives
$$
 \PP\!\left\{\max_{x\in\mathcal N_n}
 \left|\sum_{i=1}^m(Z_i(x)-\E Z_i(x))\right|>\xi\right\}
 \le 2(1+2n^{p+2})^n
 \exp\!\left(-\frac{\xi^2}{2(Cn^{-p/2}+\xi n^{-ap}/3)}\right).
$$
The right side tends to zero because $p/2>1$ and $ap>1$. The truncated
sum is $mp$-Lipschitz, and $mpn^{-p-2}=o(1)$. This proves
\eqref{trunclim} on the whole sphere.

\end{proof}

\begin{lemma}\label{gram}
Under the same assumptions, with $J(x)=\{i:|\langle U_i,x\rangle|>t\}$,
$$
 \sup_{x\in\Sph^{n-1}}\sum_{i\in J(x)}|\langle U_i,x\rangle|^2\le1+o_{\PP}(1).
 \eqid{grambound}
$$
Moreover, with probability tending to one, $|J(x)|<\lceil 2/t^2\rceil$ for every $x$.
\end{lemma}
\begin{proof}
Let $d=\lceil2/t^2\rceil$ and write $U_i=g_i/\|g_i\|_2$ with independent standard Gaussian vectors
in $\R^n$. Because $m$ is polynomial in $n$, the Gaussian norm tail estimate and the union bound
show that $\min_i\|g_i\|_2/\sqrt n\ge1-o(1)$ with probability tending to one. Put $\eta_n=(d\log
n/n)^{1/4}$, which tends to zero because $2a<1$. For a subset $J\subset\{1,\ldots,m\}$ of size
$d$, let $G_J:\ell_2^n\to\ell_2^d$ have rows $g_i$, $i\in J$. The Gaussian singular-value estimate
\cite{Ver}*{Theorem 4.6.1}, applied to the transpose, gives
$$
 \PP(\|G_J\|>\sqrt n+C\sqrt d+C\eta_n\sqrt n)\le2e^{-c_2n\eta_n^2},
$$
with absolute constants. There are at most $m^d$ such subsets, and $n\eta_n^2/(d\log n)\to\infty$.
Thus the exceptional probability is at most $2\exp(d\log m-c_2n\eta_n^2)=o(1)$. Together with
the row-length bound this gives, for all such $J$ and all unit $x$,
$$
\begin{aligned}
 \sum_{i\in J}|\langle U_i,x\rangle|^2
 &\le \frac{\|G_Jx\|_2^2}{\min_i\|g_i\|_2^2}\\
 &\le\frac{(1+C\sqrt{d/n}+C\eta_n)^2}{(1-o(1))^2}=1+o(1).
\end{aligned}
$$
Since $dt^2\ge2$, a set $J(x)$ containing $d$ indices would contradict this estimate. Every
smaller set extends to one of size $d$ for all sufficiently large $n$, which proves
\eqref{grambound}.
\end{proof}

\begin{proof}[Proof of Theorem~\ref{th:fixed}]
Write $b_i(x)=|\langle U_i,x\rangle|^2$ and $q=p/2>1$. Lemma~\ref{trunc} gives the uniform decomposition
$$
 F(x)=\alpha\mu_p+\sum_{i\in J(x)}b_i(x)^q+o_{\PP}(1).
 \eqid{decomp}
$$
Indeed the difference between the large-coordinate contribution and the excess above the
truncation is $|J(x)|t^p\le dt^p=O(t^{p-2})\to0$ on the event of Lemma~\ref{gram}.

For any finite nonnegative sequence $v_i\le r$ with $\sum_i v_i\le S$, convexity gives
$$
 \sum_i v_i^q\le \lfloor S/r\rfloor r^q+(S-r\lfloor S/r\rfloor)^q.
 \eqid{convex}
$$
To see this, append zero coordinates until the sum can be increased to $S$ while retaining the cap
$r$. For two coordinates strictly between zero and $r$, move mass from the smaller to the larger
until one is zero or the other equals $r$. Convexity of $u\mapsto u^q$ ensures the sum of powers
does not decrease. Repetition leaves at most one coordinate strictly between zero and $r$, giving
\eqref{convex}. The right-hand side is continuous in $S\ge0$ (including the integer multiples of
$r$). Applying \eqref{convex} to $b_i(x)$, $i\in J(x)$, using Lemma~\ref{gram} and \eqref{decomp},
gives the upper bound in \eqref{fixmax}.

For the lower bound fix $0<r'<r$. Let $v_1,\ldots,v_j$ be the positive entries in the vector
formed by $\lfloor1/r'\rfloor$ copies of $r'$ and the remainder $1-r'\lfloor1/r'\rfloor$. They sum
to one. Define
$$
 w=\sum_{i=1}^j\sqrt{v_i}\,U_i,\qquad x=\frac{w}{\|w\|_2}
$$
on the event $w\ne0$, and take an arbitrary deterministic unit vector otherwise. Since $j$ is
fixed and all pairwise inner products of the first $j$ vectors tend to zero in probability by
\eqref{mom}, $\|w\|_2\to1$ and $\langle U_i,x\rangle\to\sqrt{v_i}$ in probability for $i\le j$.
Conditional on the first $j$ vectors, the remaining vectors are independent and uniform on the
sphere. The cap bound in Lemma~\ref{trunc} therefore gives, for each fixed $0<\zeta<1$,
$\PP(\max_{i>j}|\langle U_i,x\rangle|>\zeta)\le
m(2e^{-n\zeta^2/4}+e^{-c_1 n})\to0$. Thus $x$ is admissible for \eqref{fixmax} with probability
tending to one, and \eqref{decomp} gives $F(x)\ge\alpha\mu_p+H_p(r')-o_{\PP}(1)$. The function
$H_p$ is continuous at every $r>0$ by its formula in \eqref{profile}. Letting $r'\uparrow r$
completes \eqref{fixmax}.

The decomposition \eqref{decomp} also bounds the minimum from below. In any fixed direction $x$,
$\E F(x)=mc_{n,p}\to\alpha\mu_p$ and $\operatorname{Var}F(x)\le mc_{n,2p}\to0$, giving the
reverse bound. Finally, if $F(X_n)\ge\sup
F-o_{\PP}(1)$, then by \eqref{fixmax} at $r=1$ and \eqref{decomp}, $\sum_{i\in
J(X_n)}b_i(X_n)^q\ge1-o_{\PP}(1)$. On the other hand,
$$
 \sum_{i\in J(X_n)}b_i(X_n)^q\le(\max_i b_i(X_n))^{q-1}\sum_{i\in J(X_n)}b_i(X_n)
 \le(\max_i b_i(X_n))^{q-1}(1+o_{\PP}(1)).
$$
Both ends tend to one in probability. Since $q>1$, this forces $\max_i|\langle U_i,X_n\rangle|\to1$ in probability.
\end{proof}

\subsection*{Optimizing the Euclidean structure}
A large value in one direction can be absorbed by a change of Euclidean structure. Large values
on a family of directions with nearly isotropic covariance constrain every comparison ellipsoid.
We will use this observation at both fixed and logarithmic exponents.

\begin{lemma}\label{ellipsoid}
Let $f$ be a norm on $\R^k$, let $\sigma_k$ be normalized surface measure on $\Sph^{k-1}$, and let
$v_1,\ldots,v_M\in\Sph^{k-1}$. Suppose that $f(v_i)\ge A>0$ for all $i$ and that, for some
$\eta\ge0$,
$$
 \frac{k}{M}\sum_{i=1}^M v_i v_i^{\top}\preceq (1+\eta)I_k.
 \eqid{covupper}
$$
Here $I_k$ is the identity matrix and $P\preceq Q$ means that $Q-P$ is positive semidefinite. Then
$$
 d_{\rm BM}((\R^k,f),\ell_2^k)^2\ge
 \frac{A^2}{(1+\eta)\int_{\Sph^{k-1}} f(x)^2\,d\sigma_k(x)}.
 \eqid{elllower}
$$
\end{lemma}
\begin{proof}
Every Euclidean norm on $\R^k$ has the form $q(x)=\sqrt{x^{\top}Qx}$ for a positive definite
symmetric matrix $Q$. After multiplying $q$ by a scalar, the distortion of the corresponding
isomorphism is the smallest $D$ such that
$$
 q(x)\le f(x)\le Dq(x)\qquad(x\in\R^k).
 \eqid{ellcompare}
$$
Integrating the first inequality in \eqref{ellcompare} after squaring, and using $\int
xx^{\top}\,d\sigma_k(x)=I_k/k$, gives
$$
 \frac{\operatorname{tr}Q}{k}\le\int_{\Sph^{k-1}}f(x)^2\,d\sigma_k(x).
$$
The second inequality in \eqref{ellcompare} gives
$$
 \frac{A^2}{D^2}\le\frac1M\sum_{i=1}^M v_i^{\top}Qv_i
 =\frac1M\operatorname{tr}\left(Q\sum_{i=1}^M v_i v_i^{\top}\right)
 \le(1+\eta)\frac{\operatorname{tr}Q}{k}.
$$
The last inequality follows from \eqref{covupper} and positivity of $Q$. Combining these estimates
and taking the infimum over the Euclidean norms proves \eqref{elllower}.
\end{proof}

\begin{corollary}\label{cor:fixed}
Let $p>2$ be fixed. If $E_N\sim\mathcal G_{N,k_N}$ and $k_N/N^{2/p}\to\beta\in(0,\infty)$, then
$$
 \Delta_p(E_N)\xrightarrow{\PP}\left(1+\frac{\beta^{p/2}}{\mu_p}\right)^{1/p}.
 \eqid{fixedlaw}
$$
The same limit holds for $\mathcal D_p(E_N)$. For $D>1$ and fixed $0<\delta<1$, the two critical dimensions satisfy
$$
 \lim_{N\to\infty}\frac{d(N,p,D,\delta)}{N^{2/p}}
 =\lim_{N\to\infty}\frac{d_{\rm iso}(N,p,D,\delta)}{N^{2/p}}
 =[\mu_p(D^p-1)]^{2/p}.
 \eqid{fixdim}
$$
\end{corollary}

\begin{proof}[Proof of Corollary~\ref{cor:fixed}]
Let $G$ be an $N\times k_N$ standard Gaussian matrix, with rows $g_i=\|g_i\|_2U_i$. It has full
column rank almost surely, and left orthogonal invariance makes its range Haar-distributed on
the Grassmannian. In these coordinates the inherited Euclidean norm is $x\mapsto\|Gx\|_2$.
Since $k_N\asymp N^{2/p}$ and $p$ is fixed, the Gaussian norm tail bound and union bound give
$\max_i|\|g_i\|_2/\sqrt{k_N}-1|\to0$ in probability. Hence, uniformly on the sphere,
$k_N^{-p/2}\|Gx\|_p^p=(1+o_{\PP}(1))\sum_i|\langle U_i,x\rangle|^p$, with multiplicative
inequalities understood on an event of probability tending to one. Apply Theorem~\ref{th:fixed}
with $n=k_N$, $m=N$, and $\alpha=\beta^{-p/2}$. The extrema of the normalized $p$th power converge
to $\beta^{-p/2}\mu_p+1$ and $\beta^{-p/2}\mu_p$. Since $k_N/N\to0$, the Gaussian singular-value
estimate gives $\|Gx\|_2/\sqrt N\to1$ uniformly in probability. Taking the ratio of the extrema
and the $p$th root gives \eqref{fixedlaw}.
For the Banach--Mazur lower bound, put $n=k_N$, $m=N$, $\alpha=\beta^{-p/2}$ and $f(x)=F(x)^{1/p}$
on $\R^n$, extending the formula for $F$ homogeneously. Almost surely the vectors span $\R^n$. By
Jensen's inequality for the concave function $t\mapsto t^{2/p}$ and \eqref{mom},
$$
 \int_{\Sph^{n-1}} f(x)^2\,d\sigma_n(x)
 \le (mc_{n,p})^{2/p}\longrightarrow(\alpha\mu_p)^{2/p}.
 \eqid{fixedmean}
$$
From \eqref{decomp}, applied at each $U_i$, we have $\min_i F(U_i)\ge\alpha\mu_p+1-o_{\PP}(1)$,
since $i\in J(U_i)$ and $|\langle U_i,U_i\rangle|=1$. Write $U_i=g_i/\|g_i\|_2$ with independent
Gaussian rows and let $A:\ell_2^n\to\ell_2^m$ be the matrix with rows $g_i$. Then
$$
 \frac{n}{m}\sum_i U_i U_i^{\top}\preceq
 \frac{n}{\min_i\|g_i\|_2^2}\frac{A^{\top}A}{m}.
$$
The scalar factor tends to one in probability by Gaussian norm concentration and the union bound,
and $\|A\|^2/m\to1$ in probability by the Gaussian singular-value bound, since $m/n\to\infty$.
Thus \eqref{covupper} holds with $\eta=o_{\PP}(1)$, chosen nonnegative. Lemma~\ref{ellipsoid} and
\eqref{fixedmean} give
$$
 d_{\rm BM}((\R^n,f),\ell_2^n)\ge
 \left(1+\frac1{\alpha\mu_p}\right)^{1/p}-o_{\PP}(1).
$$
The norm $k_N^{-1/2}\|Gx\|_p$ differs from $f(x)$ by multiplicative factors $1+o_{\PP}(1)$
uniformly in $x\ne0$, by the row length estimate at the beginning of the proof. The identity
between these two normed spaces therefore has distortion tending to one in probability. By the
multiplicative triangle inequality for the Banach--Mazur distance, the same lower bound holds for
$\mathcal D_p(E_N)$. The upper bound follows from $\mathcal D_p(E_N)\le\Delta_p(E_N)$ and
\eqref{fixedlaw}. For a fixed $D>1$, the limit in \eqref{fixedlaw} is continuous and strictly
increasing in $\beta$, and equals $D$ at $\beta=[\mu_p(D^p-1)]^{2/p}$. Taking dimensions just
below and just above this value times $N^{2/p}$, the two success probabilities tend respectively
to one and zero. Monotonicity of both distances proves \eqref{fixdim}.
\end{proof}

\subsection*{Spherical witnesses for the Euclidean identity}
The Dvoretzky--Rogers theorem rules out an absolutely $p$-summing identity in infinite dimension
for $1\le p<\infty$ \cite{DR}. In real Euclidean dimension $n$, Gordon's formula gives
$\pi_p(I_n)=c_{n,p}^{-1/p}$ \cite{Gordon}. For independent uniform unit vectors
$U_1,\ldots,U_m$, the proportion of this norm recovered by the sample is
$$
 R_{n,m,p}=\frac{m^{1/p}}{\pi_p(I_n)\sup_{x\in\Sph^{n-1}}(\sum_i|\langle U_i,x\rangle|^p)^{1/p}}.
 \eqid{ratio}
$$
\begin{corollary}\label{cor:js}
Fix $p>2$. If $m/n^{p/2}\to\alpha\in(0,\infty)$, then
$$
 R_{n,m,p}\xrightarrow{\PP}\left(\frac{\alpha\mu_p}{1+\alpha\mu_p}\right)^{1/p}.
 \eqid{jslaw}
$$
In particular, for $0<\theta<1$, the critical coefficient for recovery of a proportion $\theta$
of $\pi_p(I_n)$ by independent spherical test vectors is
$$
 \alpha_{p,\theta}=\frac{\theta^p}{\mu_p(1-\theta^p)}.
 \eqid{jsconstant}
$$
For $m/n^{p/2}\to\alpha<\alpha_{p,\theta}$, $\alpha>0$, the probability of recovery tends to
zero. For $\alpha>\alpha_{p,\theta}$ it tends to one.
\end{corollary}

\begin{proof}
Substituting Gordon's formula in \eqref{ratio} gives
$$
 R_{n,m,p}^p=\frac{mc_{n,p}}{\sup_{x\in\Sph^{n-1}}F(x)}.
$$
By \eqref{mom} and Theorem~\ref{th:fixed}, the numerator and denominator converge to $\alpha\mu_p$
and $1+\alpha\mu_p$, respectively. This gives \eqref{jslaw}. Solving
$(\alpha\mu_p/(1+\alpha\mu_p))^{1/p}=\theta$ gives \eqref{jsconstant} and the stated
strict-threshold assertions.
\end{proof}

When the range embeds isometrically in $L_p$, averaging over rotations compares an arbitrary
operator on $\ell_2^n$ with the identity on the same domain.

\begin{proposition}[Transfer to an $L_p$ range]\label{prop:gaussupper}
Let $H=\ell_2^n$, $1\le p<\infty$, and $T:H\to L_p(\Omega,\lambda)$
be a nonzero real operator. If $I_n$ denotes the Euclidean identity,
then for every $m$,
$$
 \frac{\pi_p^{(m)}(T)}{\pi_p(T)}
 \ \ge\ \frac{\pi_p^{(m)}(I_n)}{\pi_p(I_n)} .
$$
The identity attains the worst ratio over this class. For fixed $p>2$
and $m/n^{p/2}\to\alpha>0$, it follows that
$$
 \limsup_{n\to\infty}\ \sup_{T}
 \frac{\pi_p(T)}{\pi_p^{(m)}(T)}
 \le \left(1+\frac{1}{\alpha\mu_p}\right)^{1/p}.
$$
Here the supremum includes all nonzero maps with the indicated domain
and ranges that embed linearly and isometrically into a real $L_p$.
\end{proposition}

\begin{proof}
Choose representatives $a_i=Te_i$ and set $a=(a_1,\ldots,a_n)$.
Pointwise evaluation gives $Tx=\langle a,x\rangle$ almost everywhere.
Let $B=\int |a|_2^p\,d\lambda$. For every finite tuple, its $p$th
power output sum is at most $B$ times its weak denominator to the
$p$th power. Conversely the spherical probability as a family of
witnesses has numerator $c_{n,p}B$ and weak denominator $c_{n,p}$,
both in $p$th powers. Finitely supported spherical probabilities
approach these values, uniformly over unit test directions by
compactness and equicontinuity. Thus $\pi_p(T)^p=B$.
The same spherical argument gives $\pi_p(I_n)^p=1/c_{n,p}$.

Rotate any tuple $(x_j)$ by a common Haar orthogonal transformation
$Q$. Its weak denominator stays fixed, and
$$
 \E_Q\sum_j\|TQx_j\|_p^p
 =c_{n,p}B\sum_j|x_j|_2^p .
$$
Some rotation attains at least this mean. Taking the supremum proves
the comparison. The Euclidean norm has an isometric $L_p$ realization
by spherical integration, so equality holds for a member of the class.
Finally, $m$ independent spherical witnesses for $I_n$ have numerator
$m$ and weak denominator $1+\alpha\mu_p+o_{\PP}(1)$, in $p$th
powers, by Theorem~\ref{th:fixed}. This gives the asserted deterministic
limsup by taking a realization for each dimension.
\end{proof}

In particular, every $\alpha>[\mu_p((1+\varepsilon)^p-1)]^{-1}$
suffices in large dimension for this class. The lower bound of
Corollary~\ref{cor:necessary} comes from operators with polyhedral
ranges. Within the class of $L_p$ ranges, the optimal coefficient is exactly the one for the
Euclidean identity with optimally chosen witnesses.

\section{Adaptive selection from a spherical sample}\label{sec:adaptive}

How much can be gained by drawing more directions than are needed and choosing the best weighted
subfamily afterwards? The diagonal contribution in Corollary~\ref{cor:js} persists even when the
choice depends on the entire sample. For integers $1\le
m\le M$ define
$$
 \mathcal W_{m,M}=\{w\in[0,1]^M:\textstyle\sum_iw_i=1,\ |\operatorname{supp}w|\le m\}.
$$
For independent uniform points $U_1,\ldots,U_M\in\Sph^{n-1}$, put
$$
 D(w)=\sup_{x\in\Sph^{n-1}}\sum_{i=1}^M w_i|\langle U_i,x\rangle|^p,
 \qquad D_{n,m,M}=\min_{w\in\mathcal W_{m,M}}D(w).
 \eqid{eq:adaptive-def}
$$
The minimum exists by compactness. A weight $w_i$ corresponds to the test vector $w_i^{1/p}U_i$
for the Euclidean identity. For the resulting family the numerator is one and the weak denominator is $D(w)^{1/p}$.

\begin{mainresult}[Adaptive spherical witnesses]\label{th:adaptive}
Let $p>2$ and $m/n^{p/2}\to\alpha\in(0,\infty)$. If $M\ge m$ and $\log M=o(n)$, then
$$
 n^{p/2}D_{n,m,M}\xrightarrow{\PP}\mu_p+\alpha^{-1}.
 \eqid{eq:adaptive-limit}
$$
More precisely, for each $\xi>0$ there are $c,c'>0$, depending only on $p,\alpha,\xi$, such that
whenever $m\le M\le\exp(cn)$, for all sufficiently large $n$,
$$
 \PP\left(n^{p/2}D_{n,m,M}<\mu_p+\alpha^{-1}-\xi\right)
 \le \exp(-c'nm).
 \eqid{eq:adaptive-prob}
$$
\end{mainresult}

\begin{proof}
Assigning weight $1/m$ to the first $m$ vectors gives
$$
 n^{p/2}D_{n,m,M}\le \frac{n^{p/2}}m\sup_{x\in\Sph^{n-1}}\sum_{i=1}^m|\langle U_i,x\rangle|^p.
$$
By Theorem~\ref{th:fixed} with $r=1$, the right-hand side converges in probability to
$\alpha^{-1}(\alpha\mu_p+1)$. It remains to bound the probability that any weighted subfamily
improves this limit by $\xi$. For fixed $p,\alpha,\xi$, choose $C_0$ so that every $w$ with
$n^{p/2}D(w)\le\mu_p+\alpha^{-1}$ satisfies
$$
 \max_iw_i\le C_0/m,
 \qquad \left\|\sum_iw_iU_iU_i^{\top}\right\|\le C_0/n.
 \eqid{eq:ad-good}
$$
Here and below $\|A\|=\|A:\ell_2^d\to\ell_2^e\|$ for an $e\times d$ matrix $A$. The first bound
follows from $w_i\le D(w)$, obtained by taking $x=U_i$ in \eqref{eq:adaptive-def}. The second
follows from the probability-weighted H\"older inequality
$$
 \sum_iw_i|\langle U_i,x\rangle|^2
 \le\left(\sum_iw_i|\langle U_i,x\rangle|^p\right)^{2/p}.
$$
The choice of $C_0$ is possible because $m/n^{p/2}\to\alpha$.

Fix $B>0$, let $\phi_B(t)=\min\{|t|^p,B\}$ for $t\in\R$, and let $L_B$ be its Lipschitz constant.
On each support of size $m$, padded with zero weights if necessary, define
$$
 H_B(w,V)=\sum_{i\ne j}w_iw_j\phi_B(\sqrt n\langle v_i,v_j\rangle),
$$
where $V$ is the $n\times m$ matrix of column vectors $v_i\in\Sph^{n-1}$. For a fixed weight
vector $w$ with $\max_iw_i\le 2C_0/m$, write $W=\operatorname{diag}(w_i)$. Consider the set
$$
 \mathcal A_w=\{V\in(\Sph^{n-1})^m:\|VWV^{\top}\|\le C_1/n\},
$$
where $C_1\ge 2C_0$ will be chosen below. If $V,V'\in\mathcal A_w$, the inequality $\|AB\|_{\rm
HS}\le\|A\|\|B\|_{\rm HS}$ gives
$$
\begin{aligned}
 &\|W^{1/2}(V^{\top}V-V'^{\top}V')W^{1/2}\|_{\rm HS}\\
 &\quad\le (\|VW^{1/2}\|+\|V'W^{1/2}\|)\|(V-V')W^{1/2}\|_{\rm HS}\\
 &\quad\le 2\sqrt{2C_0C_1/(nm)}\,\|V-V'\|_{\rm HS}.
\end{aligned}
$$
Let $A=V^{\top}V-V'^{\top}V'$. By the Lipschitz bound for $\phi_B$ and Cauchy--Schwarz,
$$
\begin{aligned}
 |H_B(w,V)-H_B(w,V')|
 &\le L_B\sqrt n\sum_{i\ne j}w_iw_j|A_{ij}|\\
 &\le L_B\sqrt n\left(\sum_{i,j}w_iw_j|A_{ij}|^2\right)^{1/2}\\
 &= L_B\sqrt n\|W^{1/2}AW^{1/2}\|_{\rm HS},
\end{aligned}
$$
since $\sum_{i,j}w_iw_j=1$. Consequently,
$$
 |H_B(w,V)-H_B(w,V')|\le 2L_B\sqrt{2C_0C_1/m}\,\|V-V'\|_{\rm HS}.
 \eqid{eq:energy-lip}
$$

For every fixed admissible weight vector, the covariance bound defining $\mathcal A_w$ holds
with high probability. To see this, use \eqref{mom} in a fixed direction $x\in\Sph^{n-1}$ to
obtain
$$
 n^j\int |\langle U_1,x\rangle|^{2j}\,d\PP
 =\frac{n^j(2j-1)!!}{n(n+2)\cdots(n+2j-2)}\le (2j-1)!!\quad(j\ge1).
$$
Expanding the exponential in its nonnegative power series, for $0\le t\le1/4$ we obtain
$$
 \int\exp(tn|\langle U_1,x\rangle|^2)\,d\PP\le (1-2t)^{-1/2}\le e^{2t}.
$$
Taking $a=m/(8C_0)$, we have $aw_i\le1/4$. Independence gives
$$
 \int\exp\left(an\sum_iw_i|\langle U_i,x\rangle|^2\right)d\PP
 \le\prod_i e^{2aw_i}=e^{2a}.
$$
Markov's inequality therefore gives
$$
 \PP\left(n\sum_iw_i|\langle U_i,x\rangle|^2>R\right)\le e^{-(R-2)a}\quad(R>2).
$$
A $1/4$-net of the unit sphere has cardinality at most $9^n$. For a positive semidefinite matrix
$Q$, its norm is at most twice the largest value of $x^{\top}Qx$ on such a net, since
subtracting a net point from a unit vector changes the quadratic form by at most $2\|Q\|/4$. Since
$m/n\to\infty$, choosing $C_1$ sufficiently large gives
$$
 \PP(V\notin\mathcal A_w)\le e^{-c_1m}
 \eqid{eq:ad-covprob}
$$
for some $c_1>0$, uniformly in the fixed weights.

The set $\mathcal A_w$ is nonempty by \eqref{eq:ad-covprob} for large $n$. Let
$L_0=2L_B\sqrt{2C_0C_1/m}$. The function on the whole product of spheres given by
$$
 V\longmapsto \inf_{V'\in\mathcal A_w}\{H_B(w,V')+L_0\|V-V'\|_{\rm HS}\}
$$
is an $L_0$-Lipschitz extension by \eqref{eq:energy-lip}. Projection of its values onto $[0,B]$
gives an extension $\widetilde H_B(w,\cdot)$ with the same Lipschitz bound. The
spherical logarithmic Sobolev inequality and its tensorization (see \cite{Ledoux}) imply that any
$L_0$-Lipschitz real function on this product, for its product Euclidean metric, satisfies
$$
 \PP\left(|f-\textstyle\int f\,d\PP|>t\right)
 \le 2\exp(-c_2nt^2/L_0^2),
$$
where $c_2>0$ is absolute and $n\ge3$. The Euclidean Lipschitz constant is also a Lipschitz
constant for the product geodesic metric. Applying this inequality to the extension gives, for
every fixed $t>0$,
$$
 \PP\left(|\widetilde H_B-\textstyle\int\widetilde H_B\,d\PP|>t\right)
 \le 2e^{-c_{B,t}nm}.
 \eqid{eq:ad-conc}
$$
The difference of the means of $H_B$ and $\widetilde H_B$ is at most $Be^{-c_1m}$ by
\eqref{eq:ad-covprob}. Writing $\gamma_1$ for standard Gaussian measure,
$$
\begin{aligned}
 \int H_B(w,V)\,d\PP
 &=\left(1-\sum_iw_i^2\right)\int_{\Sph^{n-1}}\phi_B(\sqrt n u_1)\,d\sigma_n(u)\\
 &\longrightarrow \int_\R\phi_B(t)\,d\gamma_1(t),
\end{aligned}
$$
uniformly in the fixed weights since $\sum_iw_i^2\le 2C_0/m$. The convergence of the spherical integral follows
from the Gaussian polar representation and boundedness of $\phi_B$.

Choose $B$ so that the last Gaussian integral is greater than $\mu_p-\xi/8$. Fix $0<\rho<1/2$ with
$4B\rho<\xi/8$. Round each coordinate of an admissible weight down to a multiple of $\rho/m$, then
divide by the sum of the rounded coordinates. The resulting weight $w'$ satisfies
$$
 \|w-w'\|_1\le2\rho,\qquad 0\le w_i'\le w_i/(1-\rho).
$$
Indeed the sum of the discarded masses is at most $\rho$, and normalization of the rounded vector
changes its $\ell_1$ norm by at most $\rho$. There are at most $(1+\lfloor C_0/\rho\rfloor)^m$
such rounded weights on any support for weights satisfying \eqref{eq:ad-good}. Moreover
$\max_iw_i'\le2C_0/m$ and the actual configuration lies in $\mathcal A_{w'}$ since its weighted
covariance is at most $2C_0/n$. For any configuration $V$, expanding the difference of the
products $w_iw_j$ gives
$$
 |H_B(w,V)-H_B(w',V)|\le2B\|w-w'\|_1\le4B\rho.
$$
The number of possible supports is at most $(eM/m)^m$. Apply \eqref{eq:ad-conc} with $t=\xi/8$ to
all the rounded weights and supports. The union bound gives the exceptional probability at most
$$
 2\exp\big(m\log(eM/m)+m\log(1+\lfloor C_0/\rho\rfloor)-c_{B,\xi/8}nm\big).
$$
For $c<c_{B,\xi/8}/4$ and all sufficiently large $n$, this is at most $e^{-c'nm}$ for some $c'>0$.
Outside this exceptional event, take any $w$ with $n^{p/2}D(w)\le\mu_p+\alpha^{-1}$. Its rounded
weight $w'$ satisfies $V\in\mathcal A_{w'}$, so the extension equals $H_B$ there. The estimates
for the mean and the rounding error give
$$
 H_B(w,V)\ge\mu_p-\xi/2
$$
uniformly for every $w$ satisfying the assumed bound on $D(w)$.

Finally, testing \eqref{eq:adaptive-def} at $U_j$ and taking their weighted average gives
$$
\begin{aligned}
 n^{p/2}D(w)&\ge n^{p/2}\sum_{i,j}w_iw_j|\langle U_i,U_j\rangle|^p\\
 &\ge n^{p/2}\sum_iw_i^2+H_B(w,V)\\
 &\ge n^{p/2}/m+\mu_p-\xi/2.
\end{aligned}
$$
For sufficiently large $n$, the last expression is at least $\mu_p+\alpha^{-1}-\xi$. This proves
\eqref{eq:adaptive-prob}. The limit \eqref{eq:adaptive-limit} follows by applying it for each
fixed $\xi>0$ and using the upper bound already proved.
\end{proof}

By Gordon's formula the $p$th power of the optimal recovery proportion is $c_{n,p}/D_{n,m,M}$.
For $M\ge m$ with $\log M=o(n)$, the proportion therefore converges in probability to
$$
 \left(\frac{\alpha\mu_p}{1+\alpha\mu_p}\right)^{1/p}.
$$
Thus \eqref{jsconstant} gives the same recovery threshold under adaptive selection.

\section{Random evaluation operators and optimal witnesses}\label{sec:evaluation}

Theorem~\ref{th:adaptive} optimizes among families supported on the sampled lines. To test the
sharpness of the general operator estimate, the witnesses must be allowed to point in arbitrary
directions. For an evaluation map into $\ell_\infty^M$, each witness has a coordinate at which
its image attains its norm. This associates the witness with a sampled line. We shall show that
moving away from this line in a tangential direction cannot improve the limiting quotient.

\begin{mainresult}[Unrestricted witnesses for evaluation operators]\label{thm:operator}
Fix $2<p<\infty$. Let $u_1,\ldots,u_M$ be independent and uniform on $\Sph^{n-1}$, and set
$T_n x=(\langle u_i,x\rangle)_{i=1}^M$, from $\ell_2^n$ to
$\ell_\infty^M$. Suppose $M/n^{p/2}\to\beta$ and
$m/n^{p/2}\to\alpha$, with $0<\alpha\le\beta<\infty$ and $m\le M$.
Then, in probability,
$$
 \frac{\pi_p^{(m)}(T_n)^p}{n^{p/2}}
 \longrightarrow \frac{1}{\mu_p+\alpha^{-1}},
 \qquad
 \frac{\pi_p(T_n)^p}{n^{p/2}}
 \longrightarrow \frac{1}{\mu_p+\beta^{-1}}.
$$
The witnesses in both norms range over the entire domain and may
depend on all the rows of $T_n$.
\end{mainresult}

\begin{corollary}\label{cor:necessary}
Let $2<p<\infty$. Suppose positive integers $k_n$ have the property that
$\pi_p(T)\le K\pi_p^{(k_n)}(T)$ for every real operator with
$n$-dimensional domain, where $K>1$ is fixed. Then
$$
 \liminf_{n\to\infty}\frac{k_n}{n^{p/2}}
 \ \ge\ \frac{1}{\mu_p(K^p-1)}.
$$
In particular, the coefficient in a bound
$k_n=\lceil C_{p,\varepsilon}n^{p/2}\rceil$ with $K=1+\varepsilon$
must satisfy
$$
 C_{p,\varepsilon}\ge
 \frac{1}{\mu_p((1+\varepsilon)^p-1)}
 \ \sim\ \frac{1}{p\mu_p\,\varepsilon}
 \qquad(\varepsilon\downarrow0).
$$
\end{corollary}

\begin{proof}
If the liminf were below the proposed bound, choose a positive $\alpha$
strictly between them and a subsequence with $k_n\le\lceil\alpha n^{p/2}\rceil$.
Fix any $\beta>\alpha$. Monotonicity in the number of witnesses
and Theorem~\ref{thm:operator} imply
$\mu_p+\alpha^{-1}\le K^p(\mu_p+\beta^{-1})$ along this subsequence.
Letting $\beta\to\infty$ contradicts the choice of $\alpha$.
\end{proof}

Write $s_{ik}=\sqrt n\langle u_i,u_k\rangle$,
$\rho_{ik}=\langle u_i,u_k\rangle$, and $P_i=I-u_iu_i^{\top}$.
Thus $P_i$ is orthogonal projection onto $u_i^\perp$.

\begin{lemma}\label{lem:sphere}
In the setting of Theorem~\ref{thm:operator}, choose $1/p<a<1/2$ and
$t=n^{-a}$. For every cyclic interval $I$ of $m$ indices,
$$
 \sum_{i\in I}\min\{|\langle u_i,y\rangle|^p,t^p\}
 =\alpha\mu_p+o_{\PP}(1)
$$
uniformly in $I$ and $|y|_2=1$. The full sample satisfies the same
assertion with $\beta$ in place of $\alpha$. Moreover,
$$
 \sup_{|y|_2=1}
 \sum_{\{i:\,|\langle u_i,y\rangle|>t\}}
             |\langle u_i,y\rangle|^2\le1+o_{\PP}(1).
$$
In particular, the maximum of the untruncated sum on each interval is
$1+\alpha\mu_p+o_{\PP}(1)$, uniformly over the intervals.
\end{lemma}

\begin{proof}
For a fixed interval the proof of Lemma~\ref{trunc} gives, at each point of its $n^{-p-2}$-net, a
probability at most $2\exp(-c_{\xi,p}n^{\min(p/2,ap)})$ for a deviation by a fixed $\xi>0$. There
are $M$ intervals and the net has at most $(1+2n^{p+2})^n$ points. Since $\min(p/2,ap)>1$ and
$M=O(n^{p/2})$, the union bound still tends to zero. The Lipschitz bound $pM$ extends this
estimate to the sphere. The same argument for the full sample replaces $\alpha$ by $\beta$. Apply
Lemma~\ref{gram} to the full sample to control the squared mass above the truncation,
simultaneously on all intervals. At every unit vector the excess is at most $1+o_{\PP}(1)$ since
$|\langle u_i,y\rangle|^p\le|\langle u_i,y\rangle|^2$. At a row belonging to an interval its
diagonal term contributes an excess of $1-t^p$. This proves the last assertion uniformly.
\end{proof}

\subsection*{Off-diagonal correlations}
\begin{lemma}\label{lem:pair}
Suppose $M,n\to\infty$. For any fixed $K\ge1$ and any continuous function $h:\mathbb R\to\mathbb
R$ of polynomial growth, with $G$ a standard real Gaussian,
$$
 \sum_{i\ne k}V_iW_k h(s_{ik})=\E h(G)+o_{\PP}(1)
$$
uniformly over probability vectors $V,W$ with
$\max_i\{V_i,W_i\}\le K/M$.
\end{lemma}

\begin{proof}
Let $b_n=\E h(s_{12})$ and let $B$ be the symmetric matrix with
$B_{ik}=h(s_{ik})-b_n$ off the diagonal and zero diagonal. Then
$\E\operatorname{tr}(B^4)/M^4\to0$. Indeed, the cycles with a repeated
index contribute $O(M^{-1})$ after division by $M^4$, by
\eqref{mom} and H\"older's inequality. For four distinct
indices, the four normalized inner products around the cycle converge
jointly to independent standard Gaussians. To see this, replace the
four spherical vectors by independent Gaussian vectors divided by
their lengths. The four coordinate products associated with the four
edges are centered, have covariance matrix $I_4$, and form independent
copies as the coordinate varies. The multivariate central limit theorem
and concentration of the four lengths prove the assertion. The moment
bounds in \eqref{mom} ensure uniform integrability, so the
centered cycle expectation tends to zero and $b_n\to\E h(G)$.
Thus $\|B\|_{\op}/M=o_{\PP}(1)$. Both weight vectors have Euclidean norm
at most $\sqrt{K/M}$, and their diagonal overlap is at most $K/M$.
Applying the operator bound to $V^{\top}BW$ proves the lemma.
\end{proof}

\begin{lemma}\label{lem:tangent}
Suppose $n\to\infty$ and $M/n\to\infty$. Fix $K\ge1$. Set $g(t)=p\,t|t|^{p-2}$, with $p>2$
fixed. For the independent spherical rows above, uniformly over probability vectors $W$ with
$\max_i W_i\le K/M$,
$$
  \sum_j W_j
  \left|\sqrt n\sum_i W_i g(s_{ij})P_j u_i\right|_2^2
  =o_{\PP}(1).
$$
\end{lemma}

\begin{proof}
Since $W_j\le K/M$, it suffices to show that the unweighted mean of the squared norms tends to
zero uniformly in $W$. Expanding that mean gives
$$
 \sum_{i,k} W_iW_k\,\frac1M\sum_j
 \big(n\rho_{ik}-s_{ij}s_{kj}\big)g(s_{ij})g(s_{kj}).
$$
Terms with $j=i$ or $j=k$ vanish identically. For $i=k$ the remaining
contribution is nonnegative and is bounded by
$$
 \frac{K^2 n}{M^3}\sum_{i\ne j}g(s_{ij})^2.
$$
Its expectation is $O_{p,K}(n/M)$ by \eqref{mom}.

Fix distinct $i,k$ and condition on $u_i,u_k$. For an independent
uniform $v$, write $X=\sqrt n\langle u_i,v\rangle$,
$Y=\sqrt n\langle u_k,v\rangle$, and $\rho=\rho_{ik}$. Define
$$
 A_n(\rho)=\E_v[g(X)g(Y)],\qquad
 D_n(\rho)=\E_v[Xg(X)Yg(Y)].
$$
The two averages over $j\ne i,k$, with divisor $M$, can be replaced by
$(1-2/M)A_n(\rho)$ and $(1-2/M)D_n(\rho)$. Conditional independence over $j$, H\"older's
inequality and \eqref{mom} bound the second moment of each centered average by $C_p/M$,
uniformly in the conditioning vectors. To take the supremum over the weights, first use
$W_iW_k\le K^2/M^2$ pointwise. Cauchy--Schwarz conditional on $u_i,u_k$ then bounds the expected
total absolute error by
$$
 \frac{C_{p}K^2}{\sqrt M}\big(n\,\E|\rho_{12}|+1\big)
 =O_{p,K}\big(\sqrt{n/M}\big).
$$
This bound is valid for weights chosen after observing all the rows.

For a standard real Gaussian $G$, put $a=\E[Gg(G)]$. The cancellation occurs between the linear
Hermite coefficient of $g$ and the constant coefficient of $t\mapsto tg(t)$. For standard
Gaussians $Z,Z'$ of correlation $\rho$, the Hermite expansion gives
$$
 \E[g(Z)g(Z')]=a^2\rho+O_p(|\rho|^3),\qquad
 \E[Zg(Z)Z'g(Z')]=a^2+O_p(\rho^2).
$$
Indeed, if $f=\sum_\ell d_\ell H_\ell$ in the orthonormal Hermite basis, then
$\E[f(Z)f(Z')]=\sum_\ell d_\ell^2\rho^\ell$ by the Gaussian generating function and $L_2$
approximation. Parity gives the cubic remainder for $g$ and the quadratic remainder for $tg(t)$.
Since the two products
are homogeneous of degrees $2p-2$ and $2p$, Gaussian polar coordinates
give the same two formulas for $A_n,D_n$ multiplied respectively by
$\gamma_{n,2p-2}$ and $\gamma_{n,2p}$, where
$\gamma_{n,s}=n^{s/2}/\E|G_n|_2^s\to1$ for a standard $G_n$ in $\R^n$.

By Lemma~\ref{lem:pair} with $h(t)=t^2,t^4$, uniformly in $W$,
$$
 n\sum_{i\ne k}W_iW_k\rho_{ik}^2=1+o_{\PP}(1),\qquad
 n\sum_{i\ne k}W_iW_k\rho_{ik}^4=o_{\PP}(1).
$$
Also $\sum_{i\ne k}W_iW_k=1+O(K/M)$. Substitution now gives
$$
 \sum_{i\ne k}W_iW_k\{n\rho_{ik} A_n(\rho_{ik})-D_n(\rho_{ik})\}
 =a^2\left(n\sum_{i\ne k}W_iW_k\rho_{ik}^2-1\right)+o_{\PP}(1)
 =o_{\PP}(1).
$$
Together with the diagonal estimate and the conditional variance
bound, this proves the assertion.
\end{proof}

\begin{proof}[Proof of Theorem~\ref{thm:operator}]
Consider an arbitrary finite family $(x_j)$ with positive numerator.
Discard zero images, which can only decrease its weak denominator.
Set $r_j=\|T_nx_j\|_\infty^p/\sum_l\|T_nx_l\|_\infty^p$. Choose an
index $i(j)$ and a sign so that
$z_j=\pm x_j/\|T_nx_j\|_\infty$ satisfies
$\langle u_{i(j)},z_j\rangle=1$. The inverse $p$th power of the
quotient for the resulting family is
$$
 D=\sup_{|y|_2=1}\sum_j r_j|\langle y,z_j\rangle|^p .
$$
Write $W_i=\sum_{j:i(j)=i}r_j$. Testing $D$ at $u_i$ gives $W_i\le D$.
It suffices for the lower estimate on $D$ to consider
$D\le C_0 n^{-p/2}$, with $C_0>\mu_p+\alpha^{-1}+1$ fixed.
Then $\max_iW_i\le K/M$ for a fixed $K=K(p,\alpha,\beta,C_0)$.
Spherical averaging gives $c_{n,p}\sum_jr_j|z_j|^p\le D$. By \eqref{mom}, the assumed upper
bound for $D$ controls $\sum_j r_j|z_j|^p$ uniformly. Since $p>2$ and $\langle
u_{i(j)},z_j\rangle=1$, it follows that
$$
  \sum_j r_j|z_j-u_{i(j)}|_2^2
  =\sum_jr_j(|z_j|_2^2-1)=O(1).
$$

Average the evaluations at $y=u_k$ with weights $W_k$. The terms with $k=i(j)$ contribute
exactly $n^{p/2}\sum_iW_i^2$ to the lower bound for $n^{p/2}D$. For the remaining terms put
$i=i(j)$ and use convexity to write
$$
 |\sqrt n\langle u_k,z_j\rangle|^p
 \ge |s_{ki}|^p+g(s_{ki})\sqrt n\langle u_k,z_j-u_i\rangle .
$$
The sum of the constant terms with weights $r_jW_k$ is
$\mu_p+o_{\PP}(1)$ by Lemma~\ref{lem:pair}. Since $z_j-u_i$
is perpendicular to $u_i$, Cauchy--Schwarz bounds the absolute value
of the linear term by
$$
 \left(\sum_i W_i
 \left|\sqrt n\sum_k W_k g(s_{ki})P_i u_k\right|_2^2\right)^{1/2}
 \left(\sum_jr_j|z_j-u_{i(j)}|_2^2\right)^{1/2}=o_{\PP}(1).
$$
The term $k=i$ in the inner sum vanishes. Both error estimates are uniform in the witnesses and their number. Hence
$$
 n^{p/2}D\ \ge\ \mu_p+n^{p/2}\sum_iW_i^2-o_{\PP}(1).
$$
The last sum is at least $1/m$ for families of length at most $m$,
and at least $1/M$ for arbitrary families. These estimates bound both
summing norms from above.

Conversely, test with $u_1,\ldots,u_m$, or with all $M$ rows. Each
has image norm one. Lemma~\ref{lem:sphere} gives
$n^{p/2}D\to\mu_p+\alpha^{-1}$ and
$n^{p/2}D\to\mu_p+\beta^{-1}$, respectively. These two choices give the matching lower bounds
for the summing norms.
\end{proof}

\section{Sharp row sparsification for the spherical model}\label{rows:section}

Let $u_1,\ldots,u_M$ be independent uniform points of $\Sph^{n-1}$. Fix $p>2$ and suppose
$m/n^{p/2}\to\alpha$, $M/n^{p/2}\to\beta$, with $0<\alpha<\beta<\infty$. For $w\in[0,\infty)^M$
put
$$
 \begin{gathered}
 F_w(y)=\sum_{i=1}^M w_i|\langle u_i,y\rangle|^p,
 \qquad F=F_{\boldsymbol 1},\\
 R(w)=\sup_{|y|=1}\frac{F_w(y)}{F(y)},\qquad
 L(w)=\inf_{|y|=1}\frac{F_w(y)}{F(y)}.
 \end{gathered}
$$
The rows almost surely span $\mathbb R^n$ for all sufficiently large $n$.
Write $\mathcal W_m=\{w\ge0:|\operatorname{supp}w|\le m\}$ and define
$$
\begin{aligned}
 U_{m,M}
 &=\sup_{0\ne b\in[0,\infty)^M}
   \inf\biggl\{R(w):w\in\mathcal W_m,\
            \sum_i b_i w_i\ge\sum_i b_i\biggr\},\\
 B_{m,M}
 &=\inf_{w\in\mathcal W_m}\frac{R(w)}{L(w)} .
\end{aligned}
$$
In $U_{m,M}$ the cost $b$ is prescribed before the weights are chosen. In $B_{m,M}$ a common
rescaling of the weights leaves the quotient unchanged. We interpret it as $+\infty$ when
$L(w)=0$. Both quantities concern the $p$th powers of the norms.

\begin{theorem}\label{rows:main}
For independent uniform spherical rows and fixed $p>2$,
$$
 U_{m,M}\ \xrightarrow{\mathbb P}\
 \frac{\mu_p+\alpha^{-1}}{\mu_p+\beta^{-1}},
 \qquad
 B_{m,M}\ \xrightarrow{\mathbb P}\
 1+\frac{1}{\alpha\mu_p}.
$$
For the first limit the constant cost $b_i=1$ already gives the lower
bound. For each cyclic interval $I$ of $m$ rows, let $w^I=(M/m)\boldsymbol 1_I$.
Every such interval gives the asymptotic upper bound for $B_{m,M}$. To obtain
the upper bound for $U_{m,M}$, choose for each prescribed cost $b$ an interval
satisfying $\sum_i b_iw_i^I\ge\sum_i b_i$. Such an interval exists because
the $M$ vectors $w^I$ have average $\boldsymbol 1$.
\end{theorem}

\begin{proof}
Let $I$ range over the cyclic intervals and set
$w^I=(M/m)\boldsymbol 1_I$. For the cutoff $t$ of Lemma~\ref{lem:sphere},
write
$$
 E_I(y)=\sum_{i\in I} (|\langle u_i,y\rangle|^p-t^p)_+,
 \qquad E(y)=E_{\{1,\ldots,M\}}(y).
$$
The truncated sums over $I$ and over all rows are
$\alpha\mu_p+o_{\mathbb P}(1)$ and $\beta\mu_p+o_{\mathbb P}(1)$,
uniformly in $I$ and $|y|=1$, by Lemma~\ref{lem:sphere}. The squared-mass estimate in the same
lemma gives $0\le E_I(y)\le E(y)\le1+o_{\mathbb P}(1)$, since
$|\langle u_i,y\rangle|^p\le |\langle u_i,y\rangle|^2$.
Consequently
$$
 \frac{F_{w^I}(y)}{F(y)}
 =\frac{\beta}{\alpha}\,
   \frac{\alpha\mu_p+E_I(y)}{\beta\mu_p+E(y)}
   +o_{\mathbb P}(1).
$$
The quotient $(\alpha\mu_p+s)/(\beta\mu_p+s)$ increases in $s\ge0$.
It follows, simultaneously for all $I$, that
$$
 L(w^I)\ge \ell-o_{\mathbb P}(1),\qquad
 R(w^I)\le u+o_{\mathbb P}(1),\qquad
 \ell=\frac{\beta\mu_p}{1+\beta\mu_p},\quad
 u=\frac{\beta\mu_p+\beta/\alpha}{1+\beta\mu_p}.
 \eqid{rows:interval}
$$
Each row belongs to exactly $m$ intervals, so the average of the $w^I$
is $\boldsymbol 1$. Every cost therefore attains at least its original
value on some interval, including costs chosen after observing all the
rows. This proves the upper bounds $u$ and $u/\ell$, respectively.

For the reverse bounds normalize $\sum_iw_i=M$. In the constant-cost
problem one can decrease the weights to this normalization, and the
bilateral quotient is scale invariant. Spherical integration gives
$\int F_w=\int F$, hence $L(w)\le1\le R(w)$. By
\eqref{rows:interval} it suffices in both infima to consider
$R(w)\le C$, for a fixed $C>1+u/\ell$. The preceding estimates give
$$
 F(u_i)=1+\beta\mu_p+o_{\mathbb P}(1)
 \qquad\text{uniformly in } i,
$$
since the $i$th row contributes $1-t^p$ to $E(u_i)$.
Now $w_i\le F_w(u_i)\le C F(u_i)$, so $W_i=w_i/M$ is a probability
vector with $\max_iW_i\le C'/M$ with probability tending to one.
Write $s_{ik}=\sqrt n\langle u_i,u_k\rangle$. By Lemma~\ref{lem:pair}
with $h(s)=|s|^p$,
$$
\begin{aligned}
 R(w)\sum_i W_i F(u_i)
 &\ge \sum_i W_i F_w(u_i)\\
 &=M\sum_i W_i^2+
   M n^{-p/2}\sum_{i\ne k}W_iW_k|s_{ik}|^p
 \ge \frac{M}{m}+\beta\mu_p-o_{\mathbb P}(1).
\end{aligned}
$$
The estimates hold uniformly for the indicated weights. They give
$R(w)\ge u-o_{\mathbb P}(1)$ and establish the constant-cost lower bound.

Finally let $V$ be uniform on the indices outside
$\operatorname{supp}w$. Since $M-m$ is a fixed positive fraction of
$M$ asymptotically, $\max_i V_i\le C''/M$. Disjointness of the
supports and the bilinear assertion of Lemma~\ref{lem:pair} give
$$
 L(w)\sum_i V_iF(u_i)
 \le \sum_i V_i F_w(u_i)
 =M n^{-p/2}\sum_{i\ne k}V_iW_k|s_{ik}|^p
 =\beta\mu_p+o_{\mathbb P}(1).
$$
Thus $L(w)\le\ell+o_{\mathbb P}(1)$, completing the lower bound
for $B_{m,M}$.
\end{proof}

The smallest symmetric relative error after rescaling the selected
weights is $(B_{m,M}-1)/(B_{m,M}+1)$. Its limit is
$(1+2\alpha\mu_p)^{-1}$. The cyclic construction also preserves any
$r$ prescribed costs exactly on at most $(r+1)m$ rows, with the same
bounds \eqref{rows:interval}. Indeed, the vector of original costs is
the average of the cost vectors of all the intervals. Carath\'eodory's
theorem represents it by a convex combination of at most $r+1$ of
them. The same convex combination of the corresponding weights preserves every prescribed cost.
Both bounds pass to the combination, since their pointwise moment inequalities are linear in the
weights.

\section{Gaussian power norms at logarithmic exponents}\label{sec:gaussian}

In \cite{PVZ}*{Section~6, item~5}, Paouris, Valettas and Zinn asked whether the factor
$\log(1/\varepsilon)$ in the denominator of the dimension estimate
$$
 \frac{\varepsilon\log N}{\log(1/\varepsilon)}
$$
is necessary when $p$ is at least of logarithmic order. Their Proposition~6.2, based on
Tikhomirov's result for the cube \cite{Tik14}, establishes the corresponding upper bound for
$p>(\log N)^2$. In the proof of their Theorem~5.6, the inherited distortion is expressed through
the condition $\lambda_E\le\|x\|_p\le(1+\varepsilon)\lambda_E$ on the Euclidean unit sphere of
$E$, as in \eqref{dist}. Lytova and Tikhomirov \cite{LT}*{Theorem~A and Corollary~B}
established the change from polynomially small to logarithmically small Gaussian variance around
$p=2\log N$ and described its consequences for random Euclidean sections. We study the distortion
itself on the scale $k\asymp\log N$. In particular, for $p\ge2\log N$ the leading coefficient is
$4$ in the joint limit of Corollary~\ref{cor:dim}, both for the inherited Euclidean norm and for
the optimal Euclidean norm on the section. Related estimates for Gaussian operator norms and norms
of random matrices with independent entries appear in \cites{LS25,LS26}.

Fix $c>0$ and consider $p\ge c\log N$. Define $B:(0,\infty]\to\R$ by
$$
 B(\lambda)=\begin{cases}
 \sqrt\lambda\exp(1/\lambda-1/2),&0<\lambda<2,\\
 \sqrt 2,&2\le\lambda\le\infty.
 \end{cases}
$$
For $b>0$ let $s(b)>b$ be the unique solution of
$$
 s-b-b\log(s/b)=2,
 \eqid{s}
$$
and set $s(0)=2$. The function $s$ is continuous and strictly increasing on $[0,\infty)$, and
tends to infinity. Indeed, for $s>b$, the function $s-b-b\log(s/b)$ is strictly increasing in $s$,
strictly decreasing in $b$, and vanishes at $s=b$. Continuity at zero follows from $b\log b\to0$.

For a Gaussian matrix $G$ with rows $g_i\in\R^k$, coordinates of size $t\sqrt{\log N}$ occur
with expected multiplicity of order $N^{1-t^2/2}$ up to logarithmic factors. Taking a $p$th root
leads to the function $t\exp((1-t^2/2)/\lambda)$, where $\lambda=p/\log N$. Its maximum for
$0\le t\le\sqrt2$ is $B(\lambda)$. A single larger coordinate can dominate in an exceptional
direction. The next theorem controls both contributions simultaneously over the sphere.

\begin{mainresult}[Uniform Gaussian norm]\label{th:uniform}
Fix $c,K>0$. For $N$ sufficiently large, $1\le k\le K\log N$ and $c\log N\le p\le\infty$, let $G$
have independent standard real Gaussian entries. With probability at least $1-C(\log N)^{-2}$,
$$
 \sup_{x\in\Sph^{k-1}}\left|\frac{\|Gx\|_p}{\sqrt{\log N}}-
 \max\left\{B(p/\log N),\frac{\|Gx\|_\infty}{\sqrt{\log N}}\right\}\right|
 \le C\frac{\log\log N}{\log N}.
 \eqid{unif}
$$
The constants depend only on $c,K$.
\end{mainresult}
By homogeneity the Euclidean term in \eqref{unif} is $B(p/\log N)|x|$ for arbitrary $x\in\R^k$.

\begin{proof}
Write $L=\log N$, $\lambda=p/L$, $h=L^{-1}$ and $a_N=\log L/L$. All constants in this section
depend only on $c,K$. For $x\in\Sph^{k-1}$ and $t\ge0$ put
$$
 M(x)=\max_i\frac{|\langle g_i,x\rangle|}{\sqrt L},\qquad
 A(x,t)=\#\{i:|\langle g_i,x\rangle|\ge t\sqrt L\}.
$$
We use a deterministic $h$-net $\mathcal N$ of the Euclidean sphere of cardinality at most
$(1+2/h)^k$, obtained by the disjoint-ball volume argument. In particular,
$$
 \log|\mathcal N|\le C L\log L.
 \eqid{net}
$$
Choose $R>\sqrt2+2$ so large that $R^2/4-K\log2/2>3$.
Since $\E e^{\|g_i\|_2^2/4}=2^{k/2}$, Markov's inequality and a union bound give
$$
 \PP\!\left(\max_i\|g_i\|_2>R\sqrt L\right)
 \le N2^{k/2}e^{-R^2L/4}\le e^{-2L}.
 \eqid{rowbound}
$$

Suppose first that $p<\infty$. For each $x\in\mathcal N$ and $t\ge0$, $A(x,t)$ is binomial with
mean $v(t)=N\PP(|\gamma|\ge t\sqrt L)$, where $\gamma$ is a standard Gaussian. The elementary
Gaussian upper tail bound gives $v(t)\le2\exp(L(1-t^2/2))$. Let
$$
 b(t)=L^3\max\{1,\exp(L(1-t^2/2))\}.
$$
Writing $q_t=\PP(|\gamma|\ge t\sqrt L)$, the binomial moment generating function satisfies
$$
 \int e^{A(x,t)}\,d\PP=(1+(e-1)q_t)^N\le\exp((e-1)v(t)).
$$ Hence, for sufficiently large $N$,
$$
 \PP(A(x,t)>b(t))\le\exp((e-1)v(t)-b(t))\le e^{-L^3/2}.
$$
Using \eqref{net} and summing over the grid $t_j=jh$, $0\le j\le\lceil R/h\rceil$, we obtain,
with probability at least $1-e^{-L^2}$,
$$
 A(x,t_j)\le L^3\max\{1,e^{L(1-t_j^2/2)}\}\quad(x\in\mathcal N,\ j\le\lceil R/h\rceil).
 \eqid{counts}
$$

For the lower bound set $t_*=\min\{\sqrt\lambda,\sqrt{2-10a_N}\}$. For large $N$, $t_*$ belongs
to a compact subinterval of $(0,\infty)$ depending only on $c$. Integration of the Gaussian
density over $[t_*\sqrt L,t_*\sqrt L+L^{-1/2}]$ gives
$$
 v(t_*)\ge C^{-1} L^{-1/2} e^{L(1-t_*^2/2)}\ge C^{-1} L^{9/2}.
 \eqid{meanlow}
$$
The lower-tail binomial Chernoff estimate $\PP(A<v/2)\le e^{-v/8}$ and \eqref{net} imply that,
except on an event of probability at most $e^{-L^2}$,
$$
 A(x,t_*)\ge \tfrac12 v(t_*)\qquad(x\in\mathcal N).
 \eqid{countlow}
$$
For $0<t\le\sqrt2$ define $f_\lambda(t)=t\exp((1-t^2/2)/\lambda)$. Its maximum on this interval
is $B(\lambda)$, attained at $\min\{\sqrt\lambda,\sqrt2\}$, because
$$
 f_\lambda'(t)=\exp((1-t^2/2)/\lambda)(1-t^2/\lambda).
$$
This derivative is uniformly bounded for $\lambda\ge c$ and $0\le t\le\sqrt2$. Thus
$f_\lambda(t_*)\ge B(\lambda)-Ca_N$. From \eqref{meanlow} and \eqref{countlow},
$$
\begin{aligned}
 \frac{\|Gx\|_p}{\sqrt L}
 &\ge t_* A(x,t_*)^{1/p}\\
 &\ge f_\lambda(t_*)\exp\left(-\frac{\tfrac12\log L+\log(2C)}{p}\right)\ge B(\lambda)-Ca_N.
\end{aligned}
$$
Since $\|Gx\|_p/\sqrt L\ge M(x)$, the desired lower bound on the net follows.

For the upper bound, assume also the event in \eqref{rowbound} on which all rows have norm at
most $R\sqrt L$. Put $z_i=|\langle g_i,x\rangle|/\sqrt L$ and partition $[0,R]$ into intervals
$[t_j,t_{j+1})$, adding the right endpoint to the last interval. The contribution of each
nonempty interval to $\sum_i z_i^p$, after taking the $p$th root, is at most
$$
 \min\{t_j+h,M(x)\}A(x,t_j)^{1/p}.
 \eqid{bin}
$$
If $t_j\le\sqrt2$, \eqref{counts} bounds \eqref{bin} by
$$
 L^{3/p}(t_j+h)e^{(1-t_j^2/2)/\lambda}\le B(\lambda)+Ca_N.
$$
If $t_j>\sqrt2$, \eqref{counts} bounds \eqref{bin} by $L^{3/p}M(x)\le M(x)+Ca_N$. Put
$A_x=\max\{B(\lambda),M(x)\}$. Summing the contributions of at most $RL+2$ intervals yields
$$
 \left(\sum_i z_i^p\right)^{1/p}\le (RL+2)^{1/p}(A_x+Ca_N).
$$
Here $(RL+2)^{1/p}\le 1+Ca_N$ since $p\ge cL$, and $A_x$ is bounded on the event
\eqref{rowbound}. Together with the lower bound this proves
$$
 \left|\frac{\|Gx\|_p}{\sqrt L}-\max\{B(\lambda),M(x)\}\right|\le Ca_N\qquad (x\in\mathcal N).
 \eqid{netbound}
$$
On the row event of \eqref{rowbound} the Lipschitz constants of $M$ and $x\mapsto\|Gx\|_p/\sqrt
L$ are bounded by $R$ and $RN^{1/p}\le Re^{1/c}$, respectively. The maximum with a constant does
not increase the Lipschitz constant of $M$. For $x\in\Sph^{k-1}$ choose $x_0\in\mathcal N$ with
$\|x-x_0\|_2\le h$. The difference between the two functions at $x$ is bounded in absolute value
by
$$
 Ca_N+(Re^{1/c}+R)h\le C'a_N.
$$
Thus \eqref{netbound} holds on the whole sphere. For $p=\infty$, use \eqref{countlow} with
$t_*=\sqrt{2-10a_N}$. It gives $M(x)\ge t_*$ on $\mathcal N$ and then $M(x)\ge\sqrt2-Ca_N$ on
the sphere. This proves \eqref{unif} also in that case. The total failure probability is at most
$e^{-2L}+2e^{-L^2}$ for large $N$.
\end{proof}

\section{Extremal row lengths and the critical dimension}\label{sec:Gaussian-dist}

\begin{mainresult}[Critical distortion at logarithmic exponents]\label{th:main}
Fix $c,K>0$. There are constants $C,N_0>0$, depending only on $c,K$, such that if $N\ge N_0$,
$1\le k\le K\log N$, $c\log N\le p\le\infty$ and $E\sim\mathcal G_{N,k}$, then, with probability
at least $1-C(\log N)^{-2}$,
$$
 \left|\Delta_p(E)-\max\left\{1,\frac{\sqrt{s(k/\log N)}}{B(p/\log N)}\right\}\right|
 \le C\frac{\log\log N}{\log N}.
 \eqid{main}
$$
The same estimate holds with $\Delta_p(E)$ replaced by $\mathcal D_p(E)$, on the same event after increasing $C$.
\end{mainresult}
In particular, if $p/\log N\to\lambda\in(0,\infty]$ and $k/\log N\to b>0$, the distortion
converges in probability to $\max\{1,\sqrt{s(b)}/B(\lambda)\}$. For $0<\lambda<2$, let
$b_0(\lambda)$ be the unique positive number for which $s(b_0(\lambda))=B(\lambda)^2$. Since
$B(\lambda)>\sqrt2$, such a number exists. The limit is one precisely when $0<b\le
b_0(\lambda)$. For $\lambda\ge2$, the limit is strictly larger than one for every $b>0$.

The scale $s(k/\log N)$ is the upper extreme of $N$ independent squared Gaussian row lengths.
The error in this estimate must be uniform even when $k/\log N$ tends to zero.
\begin{lemma}\label{rowlemma}
Fix $K>0$. If $1\le k\le K L$, $L=\log N$, and $g_1,\ldots,g_N$ are independent standard
Gaussian vectors in $\R^k$, then, with probability at least $1-CL^{-2}$,
$$
 \left|\max_i\frac{\|g_i\|_2^2}{L}-s(k/L)\right|\le C\frac{\log L}{L}.
 \eqid{rowext}
$$
The constant depends only on $K$.
\end{lemma}
\begin{proof}
Put $b=k/L$ and $I_b(u)=(u-b-b\log(u/b))/2$ for $u>b$. For $Y=\|g_1\|_2^2$, its moment
generating function is $(1-2t)^{-k/2}$ for $0<t<1/2$. For $u>b$ choose $t=(1-b/u)/2\in(0,1/2)$.
The Chernoff bound reads
$$
\begin{aligned}
 \PP(Y\ge uL)&\le e^{-tuL}(1-2t)^{-k/2}\\
 &=\exp\left(-\tfrac{L}{2}(u-b)+\tfrac{k}{2}\log(u/b)\right)
 =e^{-LI_b(u)}.
\end{aligned}
$$

For the reverse tail bound use the density
$$
 f_k(y)=\frac{y^{k/2-1}e^{-y/2}}{2^{k/2}\Gamma(k/2)},\qquad y>0.
$$
For $u>b$ in a fixed compact subset of $(0,\infty)$, Stirling's formula in the form
$\log\Gamma(v)=(v-1/2)\log v-v+O(1)$ for $v\ge1/2$ gives
$$
\begin{aligned}
 \log f_k(uL)
 &=(k/2-1)\log(uL)-uL/2-(k/2)\log2\\
 &\quad -(k/2-1/2)\log(k/2)+k/2+O(1)\\
 &=-LI_b(u)-\log(uL)+\tfrac12\log(k/2)+O(1).
\end{aligned}
$$
The last two logarithms are bounded in absolute value by $C\log L$, uniformly for $1\le k\le
KL$. On $[uL,uL+1]$ the logarithmic derivative of $f_k$ is $(k/2-1)/y-1/2$, which is uniformly
bounded. Integration on this interval gives
$$
 e^{-LI_b(u)-C\log L}\le\PP(Y\ge uL)\le e^{-LI_b(u)}.
 \eqid{chitail}
$$
For $0<b\le K$, $s(b)\in[2,s(K)]$ and $s(b)-b\ge2$ by \eqref{s}. Thus $I_b'(u)=(1-b/u)/2$ is
bounded below by a positive constant depending only on $K$ for $|u-s(b)|\le1/2$. Set
$u_\pm=s(b)\pm A\log L/L$ for a sufficiently large constant $A$. Since $I_b(s(b))=1$,
\eqref{chitail} implies $N\PP(Y\ge u_+L)\le L^{-3}$ and $N\PP(Y\ge u_-L)\ge L^3$. The union
bound and independence yield
$$
 \PP(\max_i\|g_i\|_2^2>u_+L)\le L^{-3},\quad
 \PP(\max_i\|g_i\|_2^2<u_-L)\le e^{-L^3}.
$$
This proves \eqref{rowext}.
\end{proof}

\begin{proof}[Proof of the inherited-distortion estimate in Theorem~\ref{th:main}]
Let $G:\ell_2^k\to\ell_p^N$ have independent standard Gaussian entries. Put $L=\log N$,
$a_N=\log L/L$, $b=k/L$, $\lambda=p/L$, and $M(x)=\|Gx\|_\infty/\sqrt L$ for $x\in\Sph^{k-1}$.
In this proof $C$ may vary and depends only on $c,K$. Since $\sup_x M(x)=\max_i\|g_i\|_2/\sqrt
L$, Theorem~\ref{th:uniform} and Lemma~\ref{rowlemma} give
$$
 \left|\sup_{x\in\Sph^{k-1}}\frac{\|Gx\|_p}{\sqrt L}-\max\{B(\lambda),\sqrt{s(b)}\}\right|\le Ca_N
$$
with probability at least $1-CL^{-2}$. Theorem~\ref{th:uniform} gives the corresponding lower
bound $\inf_x\|Gx\|_p/\sqrt L\ge B(\lambda)-Ca_N$. For the reverse bound, choose the
deterministic first coordinate vector $e_1\in\R^k$. The Gaussian tail bound gives
$$
 \PP(M(e_1)>\sqrt2+3a_N)\le2N\exp(-L(\sqrt2+3a_N)^2/2)\le 2L^{-3}.
$$
Since $B(\lambda)\ge\sqrt2$, another application of Theorem~\ref{th:uniform} yields
$$
 \left|\inf_{x\in\Sph^{k-1}}\frac{\|Gx\|_p}{\sqrt L}-B(\lambda)\right|\le Ca_N
$$
with the same probability bound, after changing $C$.

The range of $G$ is Haar-distributed on $\mathcal G_{N,k}$. The classical Gaussian
singular-value bounds \cite{Ver}*{Theorem 4.6.1} imply that with probability at least
$1-2e^{-c_0 L^2}$,
$$
 1-C_0\frac{\sqrt{k}+L}{\sqrt N}\le\frac{\|Gx\|_2}{\sqrt N}\le1+C_0\frac{\sqrt{k}+L}{\sqrt N}\qquad(x\in\Sph^{k-1}),
 \eqid{sv}
$$
where $c_0,C_0>0$ are absolute constants. In particular the error in \eqref{sv} is $o(a_N)$
uniformly for $k\le KL$. If $u=\sup_{x\in\Sph^{k-1}}\|Gx\|_p$ and
$v=\inf_{x\in\Sph^{k-1}}\|Gx\|_p$, and $\tau=C_0(\sqrt{k}+L)/\sqrt N<1$, then \eqref{dist} and
\eqref{sv} give
$$
 \frac{1-\tau}{1+\tau}\frac{u}{v}\le\Delta_p(\operatorname{range}G)\le\frac{1+\tau}{1-\tau}\frac{u}{v}.
$$
The minimum $v/\sqrt L$ stays bounded away from zero because
$B(\lambda)\ge\sqrt2$. Both $B(\lambda)$ and $\sqrt{s(b)}$ are bounded
above in the stated range. Dividing the two extremal estimates thus gives
$$
 \frac uv=
 \frac{\max\{B(\lambda),\sqrt{s(b)}\}}{B(\lambda)}+O_{c,K}(a_N).
$$
The change of Euclidean norm costs $O(\tau)=o(a_N)$ by \eqref{sv},
which proves \eqref{main}.

\end{proof}

\subsection*{The optimal comparison ellipsoid}

\begin{proof}[Proof of the Banach--Mazur estimate in Theorem~\ref{th:main}]
Let $G$ be the Gaussian matrix of Theorem~\ref{th:uniform}. Put $L=\log N$, $a_N=\log L/L$,
$b=k/L$, $\lambda=p/L$, $M(x)=\|Gx\|_\infty/\sqrt L$ and $f(x)=\|Gx\|_p/\sqrt L$ for $x\in\R^k$.
The constants may depend only on $c,K$. Almost surely $G$ has full column rank, so $f$ is a norm
and $d_{\rm BM}((\R^k,f),\ell_2^k)=\mathcal D_p(\operatorname{range}G)$.
For any deterministic unit vector $x$,
$$
 \PP(M(x)>\sqrt2+8a_N)\le2N\exp\big(-L(\sqrt2+8a_N)^2/2\big)\le2L^{-10}.
$$
By Tonelli's theorem and Markov's inequality, with probability at least $1-2L^{-8}$,
$$
 \sigma_k\{x:M(x)>\sqrt2+8a_N\}\le L^{-2}.
 \eqid{badmeasure}
$$
Intersect this event with those in \eqref{rowbound} and \eqref{unif}.
On the exceptional directions in \eqref{badmeasure}, $f$ is bounded
by a constant depending only on $c,K$. On the complement,
$f\le B(\lambda)+Ca_N$. Since $B(\lambda)$ is uniformly bounded,
integration gives
$$
 \int_{\Sph^{k-1}} f(x)^2\,d\sigma_k(x)\le B(\lambda)^2+Ca_N .
 \eqid{meanf}
$$

A large value of $f$ in one direction sufficed for the inherited
distortion. For the Banach--Mazur distance we need enough directions
of large norm to constrain an arbitrary comparison ellipsoid.

The lower tail estimate \eqref{chitail} and the uniform lower bound for $I_b'$ in the proof of
Lemma~\ref{rowlemma} imply that, for a sufficiently large constant $A_0$ depending only on $K$,
$$
 N\PP(\|g_1\|_2^2/L\ge s(b)-A_0a_N)\ge L^{10}.
$$
Indeed $I_b(s(b)-A_0a_N)\le1-c_K A_0a_N$, and $A_0$ can be chosen so that $c_KA_0-C\ge10$. The
number of rows above the threshold $L(s(b)-A_0a_N)$ in squared length is therefore binomial with
mean at least $L^{10}$. The binomial lower tail bound shows that at least $M_0=\lceil L^8\rceil$
rows satisfy it, except on an event of probability at most $\exp(-L^{10}/8)$ for sufficiently
large $N$. Choose the first $M_0$ such rows in their original order and denote their normalized
directions by $v_1,\ldots,v_{M_0}$. Conditional on all row lengths, these directions are
independent and uniform on $\Sph^{k-1}$, by independence of the Gaussian radial and angular
variables.

Put $H=(k/M_0)\sum_i v_i v_i^{\top}$. Using independence, $\int v_i v_i^{\top}\,d\PP=I_k/k$ and
$\operatorname{tr}((v_i v_i^{\top})^2)=1$, we have, conditionally on the lengths,
$$
 \int\|H-I_k\|_{\rm HS}^2\,d\PP=\frac{k^2-k}{M_0}.
$$
Here $\|A\|_{\rm HS}^2=\operatorname{tr}(A^{\top}A)$. Since the operator norm is at most the
Hilbert--Schmidt norm, Markov's inequality gives $H\preceq(1+L^{-1})I_k$ except with conditional
probability at most $L^2(k^2-k)/M_0\le K^2L^{-4}$. For each selected row $g_j$, its direction
$v_i=g_j/\|g_j\|_2$ satisfies
$$
 f(v_i)^2\ge |\langle g_j,v_i\rangle|^2/L\ge s(b)-A_0a_N.
$$
Applying Lemma~\ref{ellipsoid} with \eqref{meanf} gives
$$
 \mathcal D_p(\operatorname{range}G)^2\ge
 \frac{s(b)-A_0a_N}{(1+L^{-1})(B(\lambda)^2+Ca_N)}.
 \eqid{bmloglower}
$$
Since $B(\lambda)\ge\sqrt2$, $s(b)\in[2,s(K)]$, and $B$ is bounded on $[c,\infty]$, the square
root of the right-hand side of \eqref{bmloglower} is at least $\sqrt{s(b)}/B(\lambda)-Ca_N$. The
distance is at least one and at most $\Delta_p(\operatorname{range}G)$. The estimate in
\eqref{main} for the inherited Euclidean norm gives the matching upper bound. The total
probability of the exceptional events is at most $CL^{-2}$, with $C$ depending only on $c,K$.
\end{proof}

\subsection*{Inverting the distortion law}

\begin{corollary}\label{cor:dim}
Suppose $p_N/\log N\to\lambda\in(0,\infty]$. For every $D>1$ and $0<\delta<1$,
$$
 \frac{d(N,p_N,D,\delta)}{\log N}\longrightarrow b_{\lambda,D},
 \eqid{dimlim}
$$
where $b_{\lambda,D}\in(0,D^2B(\lambda)^2)$ is the unique solution of
$$
 b_{\lambda,D}\left(1+\log\frac{D^2B(\lambda)^2}{b_{\lambda,D}}\right)
 =D^2B(\lambda)^2-2.
 \eqid{beq}
$$
Independently of the existence of a limit for $p_N/\log N$, suppose $p_N\ge2\log N$, $\eps_N\downarrow0$ and
$$
 \frac{\eps_N\log N}{\log\log N}\longrightarrow\infty.
 \eqid{joint}
$$
Then, for every fixed $0<\delta<1$,
$$
 \frac{d(N,p_N,1+\eps_N,\delta)\log(1/\eps_N)}{\eps_N\log N}\longrightarrow4.
 \eqid{four}
$$
For fixed $\lambda\ge2$, the same constant $4$ is obtained in the iterated limit, first
$N\to\infty$ and then $D=1+\eps\downarrow1$. All these statements also hold for $d_{\rm iso}$ in
place of $d$.
\end{corollary}
The order $\eps\log N/\log(1/\eps)$ is classical for the cube \cite{Tik14}. Formula \eqref{four}
identifies its leading coefficient in the normalization \eqref{dist} and establishes the
necessity of the logarithmic factor for finite exponents $p\ge2\log N$ in the joint regime
\eqref{joint}. Equation \eqref{beq} also describes the transition below $2\log N$. The same
comparison applies after optimization over all Euclidean structures on the section.

\begin{proof}
Fix $D>1$. Since $s$ increases continuously from $2$ to infinity, it takes the value
$D^2B(\lambda)^2>2$ at a unique $b>0$. Substitution of $s=D^2B(\lambda)^2$ into \eqref{s} gives
\eqref{beq}. For any $b_-<b_{\lambda,D}<b_+$ with $b_->0$, Theorem~\ref{th:main}, with
$k=\lfloor b_-L\rfloor$ and $k=\lceil b_+L\rceil$, gives probabilities of distortion at most $D$
tending to one and zero, respectively. Monotonicity in $k$ and then $b_-,b_+\to b_{\lambda,D}$
give \eqref{dimlim}.

For $\lambda\ge2$ put $b_\eps=b_{\lambda,1+\eps}$, which is independent of $\lambda$. By \eqref{beq},
$$
 b_\eps\left(1+\log\frac{2(1+\eps)^2}{b_\eps}\right)=4\eps+2\eps^2.
$$
As $\eps\downarrow0$, continuity of $s^{-1}$ at $2$ yields $b_\eps\downarrow0$, and the
preceding identity reduces to $b_\eps\log(1/b_\eps)\sim4\eps$. Taking logarithms gives
$$
 \log(1/\eps)=\log(1/b_\eps)-\log\log(1/b_\eps)+\log4+o(1).
$$
Dividing by $\log(1/b_\eps)$ proves that $\log(1/\eps)/\log(1/b_\eps)\to1$. Hence
$$
 b_\eps\sim\frac{4\eps}{\log(1/\eps)}.
 \eqid{basym}
$$
This proves the iterated limit.

For the joint limit define $a_N=\log L/L$ and fix $0<\eta<1$. Put $D_\pm=1+(1\pm\eta)\eps_N$ and
$k_-=\lfloor Lb_{(1-\eta)\eps_N}\rfloor$, $k_+=\lceil Lb_{(1+\eta)\eps_N}\rceil$. Condition
\eqref{joint} and \eqref{basym} imply $k_\pm\to\infty$, since eventually $\eps_N\ge a_N$ and
$\log(1/\eps_N)\le\log L$. Moreover $a_N=o(\eps_N)$, so the additive error in
Theorem~\ref{th:main} is eventually smaller than $\eta\eps_N$. Theorem~\ref{th:main}, with $K=1$
and $c=2$, applies to these integers. By monotonicity of $s$, the deterministic term in
\eqref{main} is at most $D_-$ at $k_-$ and at least $D_+$ at $k_+$. The probabilities of
distortion at most $1+\eps_N$ thus tend to one and zero, respectively. In particular, for large
$N$,
$$
 k_-\le d(N,p_N,1+\eps_N,\delta)<k_+.
$$
Divide by $\eps_N L/\log(1/\eps_N)$, use \eqref{basym}, and let $\eta\downarrow0$ to obtain
\eqref{four}. The same argument applies to $d_{\rm iso}$ by the estimate for $\mathcal D_p$ and
its monotonicity under inclusion.
\end{proof}

\appendix

\section{Two refinements in the accuracy parameter}\label{sec:refinements}

Two additional structures give sharper bounds in the accuracy parameter. For even $p$ the moment
is a quadratic form on symmetric tensors. On a Hilbert domain the weak denominator is already a
supremum over a Euclidean ball, so the correlations of the original rows can be used directly.

\subsection*{Even exponents with prescribed costs}
\begin{proposition}\label{disc:even}
Let $p=2k$ with $k\ge1$ an integer. For $a_i\in\mathbb R^n$ and
$b_i\in\mathbb R^s$, the prescribed cost identities and moment inequalities of
Theorem~\ref{full:rows} hold for every $0<\eta<1$
with support at most
$$
 C\eta^{-2}\binom{n+k-1}{k}+128(s+1),
$$
where $C$ is absolute.
\end{proposition}

\begin{proof}
In the Euclidean symmetric $k$-tensor space put $v_i=a_i^{\otimes k}$,
$A_i=v_iv_i^{\top}$, and $B=\sum_i A_i$. We argue on the span of
the $v_i$, whose dimension $d$ is at most $\binom{n+k-1}{k}$.
At any stage of the iteration, absorb each positive weight in its matrix $A_i$ and its cost
$b_i$, let $M$ count the surviving indices, and write $B$ for the current sum of their matrices.
This convention also covers zero matrices with nonzero costs. When $d=0$, conic Carath\'eodory
gives the assertion.
When $M\ge C_0d$, take $\delta=(C_0d/M)^{1/2}$, increasing the
stopping threshold below so that $\delta\le1/2$ at every step. Define
$$
 K=\{z\in[-1,1]^M:\ \textstyle\sum_i z_i A_i\preceq\delta B\},
 \qquad Q=K\cap(-K).
$$
Since each $A_i$ is positive semidefinite, $[-\delta,\delta]^M\subset K$. For every set $S$ of
coordinate indices, independent uniform coefficients on $[-1,1]^S$ satisfy
$\sum_{i\in S} z_i A_i\preceq0$ with probability at least
$\exp(-2d)$. This is Corollary~30 of Reis and Rothvoss \cite{RR},
applied after reversing all signs. Thus every coordinate section of $K$ has at least an
$e^{-2d}$ fraction of the volume of its cube. Their Theorem~22 gives
$\operatorname{vol}_M(Q)\ge2^{-5M}$, since $M\delta^2=C_0d\ge\log_2(e^{2d})$ when
$C_0\ge2/\log2$.

Since $Q\subset[-1,1]^M$, its standard Gaussian measure is at least
$(2^{-5}/\sqrt{2\pi e})^M$. Let
$H=\{z:\sum_i z_i b_i=0\}$ for the current costs.
When $M\ge4s$ this subspace has dimension at least $3M/4$.
The constrained partial coloring theorem of Reis and Rothvoss
\cite{RRbal}*{Theorem~6 of the preprint}, with shift zero, $\beta=3/4$ and $\gamma=1/4$,
gives an absolute $C_1$ and $z\in C_1 Q\cap H\cap[-1,1]^M$
with at least $M/2$ coordinates of modulus one.
Choosing its sign makes at least $M/4$ of the weights $1+z_i$
vanish. The update preserves all costs and places the new matrix
sum between $(1-C_1\delta)B$ and $(1+C_1\delta)B$.

Iterate until the support is at most
$K_*=\lceil C\eta^{-2}\binom{n+k-1}{k}+4(s+1)\rceil$.
The support sizes immediately preceding a reduction, read in reverse
order, are at least $K_*(4/3)^j$, $j=0,1,\ldots$. The sum of their
errors $C_1\delta$ is therefore bounded by
$C_1\sqrt{C_0\binom{n+k-1}{k}/K_*}/(1-\sqrt{3/4})$.
Choose $C$ so this sum is at most $\eta/3$. The products of the
successive upper and lower factors then belong to
$[1-\eta,1+\eta]$. Testing the resulting matrix inequalities on
$x^{\otimes k}$ proves the moment estimate. The larger displayed
allowance for the cost dimension absorbs rounding.
\end{proof}

\subsection*{A Hilbert domain}

\begin{theorem}[Hilbert domains]\label{disc:hilbert}
Let $p>2$. Suppose $H$ is an $n$-dimensional Hilbert space, $n\ge1$, and $T:H\to Y$ is a linear
map into a Banach space over the same scalar field, real or complex. If
$0<\varepsilon<1$, then
$$
 \pi_p(T)\le(1+\varepsilon)\pi_p^{(m)}(T),\qquad
 m=\left\lceil C_p\varepsilon^{-2}\log(2/\varepsilon)n^{p/2}\right\rceil .
$$
\end{theorem}

\begin{proof}
We give the real argument first. Let $a_1,\ldots,a_M$ lie in the Euclidean unit ball, put
$F(y)=\sum_i|\langle a_i,y\rangle|^p$, and set $W=\sup_{|y|\le1}F(y)$. Suppose $W\ge1$ and
$M\ge64$.
Given any $b\in\mathbb R^M$, we shall find $\chi\in[-1,1]^M$
orthogonal to $b$, with at least $M/2$ saturated coordinates, such that
$$
 \sup_{|y|\le1}\left|\sum_i\chi_i|\langle a_i,y\rangle|^p\right|
       \le C_p\sqrt{W\log(2+W)} .
 \eqid{hilb:step}
$$

At a dyadic level $t=2^{-\ell}$ use the function $f_t$ of Section~\ref{sec:support-count}. Its
support is in $t/2\le|z|\le2t$ and its Lipschitz constant is at most $C_pt^{p-1}$.
Set
$$
 e=t/32,\qquad \tau=(32\sqrt W)^{-1},\qquad
 \theta=(\tau t)^2,\qquad
 \delta=\frac{\theta t^p}{4^p W}.
$$
Join distinct rows if $|\langle a_i,a_j\rangle|\ge\delta$.
The maximum closed-neighborhood size $D$ is at most $1+W\delta^{-p}$
by testing $F$ at each row. In particular
$\log(2D/e)\le C_p\log((2+W)/t)$.

Consider ordered tuples $A$ of $1\le k\le e^{-2}$ occurrences,
allowing repetitions, whose indices are connected in the third power
of this graph. Each occurrence $h$ bears a coefficient
$c_h\in[-1,1]\cap e^3\mathbb Z$. On the union $U_A$ of the closed
neighborhoods of the indices in $A$, set
$$
 s_i=\sum_{h\in A}c_h\langle a_h,a_i\rangle,\qquad
 J_A=\{i\in U_A:t/2-3e\le |s_i|\le 2t+3e\},
$$
where $a_h$ denotes the row at that occurrence. For every nonempty
$J_A$ and each dyadic radius $t/8\le R\le1$, use the seminorm
$$
 Z_{A,R}(x)=\sup_{|v|\le R}
              \left|\sum_{i\in J_A}x_i f_t(\langle a_i,v\rangle)\right|.
$$
Write $L_t=\operatorname{Lip}(f_t)\le C_pt^{p-1}$. For a standard Gaussian vector $g$, the
process in the definition of $Z_{A,R}(g)$ has squared increments bounded by those of $v\mapsto
L_t\sum_{i\in J_A}g_i\langle a_i,v\rangle$. Since $v=0$ is admissible, symmetry bounds the
expected absolute supremum by twice the expected supremum. Sudakov--Fernique comparison
therefore gives
$$
 \E Z_{A,R}(g)
 \le 2L_t R\,\E\left|\sum_{i\in J_A}g_i a_i\right|
 \le C_pt^{p-1}R\sqrt{|J_A|}.
$$
The seminorm itself is $2^pt^p\sqrt{|J_A|}$-Lipschitz on $\R^M$. We impose the symmetric convex constraints
$$
 Z_{A,R}(x)\le C_p\sqrt{|J_A|}
       \left(t^{p-1}R+t^p\sqrt{k\log(2D/e)}\right).
 \eqid{hilb:gauss}
$$
There are at most
$M(4+\ell)[k^2 D^3(1+2e^{-3})]^k$ entries at $(\ell,k)$.
This follows by rooting a spanning tree on the labeled occurrences.
Its root admits $M$ choices and each further vertex at most $D^3$.
The remaining factors account for the tree, coefficients and radius.
Gaussian concentration, with the constant in \eqref{hilb:gauss}
sufficiently large, bounds the total failure probability at
$(\ell,k)$ by $M(e/2)^{4k}$, since $k\le e^{-2}$.
Each failure probability is at most $1/2$, and their sum is less than $M/64$. The constraint
sets are symmetric and convex. By Royen's Gaussian correlation inequality \cite{Royen}, their
intersection has measure at least the product of their measures, hence at least $\exp(-M/32)$
since $\log(1-x)\ge-2x$ on $[0,1/2]$. Continuity from above includes all levels. One more
application to the cube gives a lower bound $\exp(-CM)$ for the measure of the intersection with
$[-1,1]^M$. Apply the constrained partial coloring theorem of Reis and Rothvoss
\cite{RRbal}*{Theorem~6 of the preprint} with shift zero, $\beta=3/4$ and $\gamma=1/4$. Since
$\dim b^\perp\ge3M/4$, it gives a coloring orthogonal to $b$ with at least $M/2$ saturated
coordinates. Each seminorm bound in \eqref{hilb:gauss} holds up to an additional absolute
factor.

Fix $|y|\le1$ and a level $t$. Starting with $r=y$, whenever
$|c|=|\langle r,a_i\rangle|>e$ record $i$ and
$c=\langle r,a_i\rangle$ and replace $r$ by $r-ca_i$.
Each update decreases $|r|^2$ by at least $c^2$, so there are at
most $e^{-2}$ occurrences and $\sum|c|\le e^{-1}$. Round towards
zero on the mesh $e^3$, and partition the list into components
$A$ in the third power of the graph, joining repeated indices.
Their sets $U_A$ are disjoint. At $i\in U_A$ the residual and the
rounding each cost at most $e$, and the omitted components contribute
at most $\delta/e\le e$. Outside all $U_A$ the same bounds give
$|\langle a_i,y\rangle|\le3e$. Writing $z_i=\langle a_i,y\rangle$,
we thus have
$$
 \sum_i\chi_i f_t(z_i)=\sum_A\sum_{i\in J_A}\chi_i f_t(z_i),
 \qquad
 i\in J_A\ \Longrightarrow\ t/4\le |z_i|\le3t .
$$
In particular $S=\sum_A|J_A|\le4^pt^{-p}F(y)$. We now discard
empty sets $J_A$.

Let $A_C$ denote the matrix with rows $a_i$, $i\in J_C$, for a
component $C$. Let $P_C$ project orthogonally onto the eigenspaces
of $A_C^{\top} A_C$ with eigenvalues at least $\theta$. The bound
$a_i a_i^{\top}\le A_C^{\top} A_C$ gives $|(I-P_C)a_i|\le\sqrt\theta$.
If $C\ne C'$, no edge joins $J_C$ to $J_{C'}$, since such an edge
would give a path of length at most three between their anchors.
Hence $\|A_CA_{C'}^{\top}\|\le\delta\sqrt{|J_C||J_{C'}|}$.
Every vector $v_C$ in the range of $P_C$ is $A_C^{\top} d_C$ with
$|d_C|\le|v_C|/\sqrt\theta$, by singular value decomposition.
Consequently
$$
 \left|\sum_Cv_C\right|^2
 \le\sum_C|v_C|^2+\frac{\delta}{\theta}
                \left(\sum_C \sqrt{|J_C|}\,|v_C|\right)^2
 \le (1+\delta S/\theta)\sum_C|v_C|^2
 \le2\sum_C|v_C|^2 .
$$
The adjoint inequality gives $\sum_C|P_Cy|^2\le2$.
For $i\in J_C$ the projected evaluation differs from $z_i$ by at most
$\sqrt\theta=\tau t$. Thus $|P_Cy|\ge t/4-\tau t>t/8$.
Choose the smallest dyadic $R_C\ge|P_Cy|$ in $[t/8,1]$.
Then $\sum_C R_C^2\le8$.

Evaluate \eqref{hilb:gauss} at $v=P_Cy$ on each set.
If $k_C$ is the number of occurrences in $C$, the two terms in \eqref{hilb:gauss} sum to at most
$$
 C_p t^{p-1}\sqrt{S\sum_C R_C^2}
 +C_p t^p\sqrt{S\log(2D/e)\sum_C k_C}
 \le C_p\sqrt S\big(t^{p-1}+t^pe^{-1}\sqrt{\log(2D/e)}\big).
$$
Here $\sum_C R_C^2\le8$ and $\sum_C k_C\le e^{-2}$. Replacing the projected evaluations by $z_i$
costs at most $C_p\tau t^p$ per selected index, and every such index satisfies
$t/4\le|z_i|\le3t$. This comparison
is confined to the selected sets, even if a projected evaluation
crosses an endpoint of the band. Each nonzero $|z_i|$ occurs in only
a bounded number of these enlarged dyadic bands, so the errors sum
to at most $C_p\tau F(y)$. The remaining contributions sum to at most
$$
 C_p\sqrt{F(y)}\sum_{\ell\ge0}t^{p/2-1}
      \sqrt{\log((2+W)/t)}
 \le C_p\sqrt{F(y)\log(2+W)} .
$$
The convergence uses $p>2$. Our choice of $\tau$ gives \eqref{hilb:step}.

Finally take any witness $x_1,\ldots,x_N$ of positive numerator and
put $h=8\sum_i|x_i|^p/N$. Split $x_i$ into
$L_i=\max(1,\lceil |x_i|^p/h\rceil)$ copies of $x_i/L_i^{1/p}$.
After common division by $h^{1/p}$ these are unit-ball rows,
$N\le M\le9N/8$, and their total $p$th power of lengths is $N/8$.
Spherical integration and Jensen's inequality give
$W\ge u=N/(8n^{p/2})$. Give each copy of $x_i$ the cost
$\|Tx_i\|^p/L_i$. Choose the sign of $\chi$ so that at least $M/4$
of the weights $1+\chi_i$ vanish. Absorbing the $p$th roots of the
weights into the copies and restoring the scale preserves the
numerator exactly. For $u\ge1$ the weak denominator to the power $p$
increases by a factor at most $1+C_p\sqrt{\log(2+u)/u}$, by
\eqref{hilb:step} and monotonicity of $\log(2+u)/u$.
At most $27N/32$ copies survive.

Repeat while $N>m$. Read backwards, the pre-reduction values of
$u$ are at least $u_*(32/27)^j$, $j=0,1,\ldots$, where
$u_*=m/(8n^{p/2})$. Summing the relative increases gives at most
$C_p\sqrt{\log(2+u_*)/u_*}$ by the convergence of
$\sum_{j\ge0}(1+j)^{1/2}(27/32)^{j/2}$. With the constant in $m$
sufficiently large, the product of the factors is at most
$(1+\varepsilon)^p$. Taking the supremum over initial witnesses
proves the real assertion.

For complex $H$ and $Y$, regard both spaces as real.
Write $T_{\R}$ for realification and set
$d_p=(2\pi)^{-1}\int_0^{2\pi}|\cos\theta|^p\,d\theta$. Every real
functional of norm at most one is the real part of a complex
functional of norm at most one. For a fixed family, maximize its complex weak moment over the
dual unit ball and average the real part of the maximizing functional over the phases. The real
weak denominator raised to the $p$th power is at least $d_p$ times the
$p$th power of the complex weak denominator. Hence $\pi_p^{(m)}(T_{\R})\le
d_p^{-1/p}\pi_p^{(m)}(T)$. For the unrestricted norms equality holds. Indeed, replace each $x_j$
in a complex witness family by $Q^{-1/p}e^{2\pi i k/Q}x_j$,
$1\le k\le Q$. Its numerator stays fixed, and the $p$th power of its
real weak denominator converges to $d_p$ times the $p$th power of
its complex weak denominator before this replacement. The convergence
is uniform in the functional by uniform convergence of the phase
Riemann sums for translates of $|\cos|^p$. Applying the real result to $T_{\R}$ in dimension
$2n$ and canceling $d_p^{-1/p}$ proves the complex case, with $2^{p/2}$ absorbed in $C_p$.
\end{proof}

\section*{Acknowledgments}
D. M. Pellegrino acknowledges support from CNPq through Grants No. 406457/2023-9 (CNPq/MCTI No. 10/2023), 403964/2024-5 (MCTI/CNPq No. 16/2024), and 305807/2025-0. E. V. Teixeira gratefully acknowledges support from the Grayce B. Kerr Chair funds at Oklahoma State University. This research was conducted in part under the DARPA ExpMath project ``A Human-Centered Framework for AI-Mathematician Collaboration in Research-Level Mathematics'' (Agreement No. HR0011262E029), in which E. V. Teixeira serves as a co-principal investigator. He gratefully acknowledges DARPA's support and the collaborative environment fostered by the project. The views, opinions, and findings expressed here are those of the authors and should not be interpreted as representing the official views or policies of the Department of Defense or the U.S. Government.

\section*{Declarations}
\noindent\textit{Data availability.} The paper uses no external datasets.

\smallskip
\noindent\textit{Competing interests.} The authors declare no competing interests.

\smallskip
\noindent\textit{Use of artificial intelligence.} 
During the elaboration of this article, we used ChatGPT 6.2 Astra (OpenAI) to
explore approaches and to test proposed estimates against
extremal examples. One instance was
the counting argument in Lemma~\ref{det:count}. The graphs at successive
stages may have increasingly large degrees, while the lists of selected
columns become shorter. ChatGPT helped us work out how to count each
choice using the degree at the stage where it occurs. This preserves
the gain from the shorter lists in the discrepancy estimate for
$2<p<3$. We also used ChatGPT in drafting and revising the exposition.
We treated all AI-generated suggestions as provisional, checked the
arguments retained in the paper, and take full responsibility for the
final manuscript.

\end{document}